\documentclass[11pt]{article}
\usepackage[T1]{fontenc}
\usepackage[sc]{mathpazo}
\usepackage{amsmath}
\usepackage{amsmath}		% for matrix environment
\usepackage[includefoot]{geometry}
\usepackage{amsmath}
\usepackage{amssymb}
\usepackage{amsthm}
\usepackage{array}
\usepackage{xcolor}
\usepackage{graphicx}

\usepackage{tikz}
\usepackage{tikz-3dplot}
\usepackage{pgfplots}
\usepackage{quiver}
\usepackage{tkz-graph}
\usetikzlibrary{positioning,arrows.meta,backgrounds}
\usetikzlibrary{arrows}
\usetikzlibrary{arrows.meta}
\usetikzlibrary{automata}
\usetikzlibrary{decorations.markings}
\usetikzlibrary{patterns}
\usetikzlibrary{positioning}
\usetikzlibrary{shapes.geometric}
\usepackage{pdfpages}

\usepackage{pifont}
\usepackage{forest}

\usepackage{multirow,multicol,enumitem}

\usepackage{ytableau}
\usepackage{bm}
\usepackage{bbm}

\usepackage{caption}
\usepackage{mathrsfs}
\usepackage{standalone} 
\usepackage{xstring}
\usepackage{rotating} 
\theoremstyle{plain}
\newtheorem{theorem}{Theorem}

\newtheorem{question}[theorem]{Question}
\theoremstyle{definition}
\newtheorem{definition}[theorem]{Definition}
\newtheorem{example}[theorem]{Example}
\theoremstyle{remark}

\usepackage{titlesec}

\newcommand{\name}[1]{\noindent {\it #1}\medskip}

\begin{document}

\section*{Inverse limits of various posets}

  \name{Amrita Acharyya}\\

\begin{abstract}
	It is known 
when we call a poset P, a
$\mathcal{P}$-chain permutational poset,
given a subset of permutations $\mathcal{P}$
of the symmetric group $S_{n}$. In this work, we use the same idea to study subsets of
words of length $n$, that are not necessarily
permutations, for example: especially when they are certain classes of restricted growth
functions induced by set partitions in standard form over $[n]=\{1,2\cdots n\}$. Varying $n$ only, and also varying $n$ and $k$ (the number of blocks of the set partitions) simultaneously, we can show that those posets form a projective system of trees and lattices (after giving a lattice structure in a natural way). These poset structures can be extended over signed restricted growth functions for standard type B set partitions over $\langle n\rangle=\{-1,-2,\cdots n,0,1,2\cdots n\}$ as well.
We investigate properties of the tree and lattice structures of these projective systems. In this scenario we further bring up some other posets like $\mathcal{P}$-Partition posets of snake graph of continued fractions, Ascent lattices on Dyck Paths, certain type of lattice induced by generalisec fibonnaci number and Stanley order, lattices induced by non-crossing set partitions.
\\

\end{abstract}
\section{Introduction}
In the preliminary section, we initially provide the background for set partitions of $[n]$, bijection with restricted growth functions, where the number of blocks $k$ of the set partition of $[n]$ turns out to be the maximal letter in the restricted growth function. We describe the connection of pattern avoidence classes with restricted growth functions. There has been quite a lot of research in the set partitions of $[n]$ as in \cite{MR2445243},
\cite{MR1544457}  and it has been extended to type B set partitions over $\langle n \rangle$ as in as well. See \cite{MR1644458}  
\cite{MR1644469} 
\cite{MR1644470} 
\cite{MR1644471}
. The idea of restricted growth functions are extended in case of of $\langle n \rangle$ utilizing type B set partitions. We also, provide some initial background of Projective Fraisse families of graphs and tress with confluent and monotone bonding epimorphisms.  We notice that the idea of chain permutational posets studied by Rodica Simion and Frank W Schmidt in \cite  {MR1644460} can be extended utilizing words of length $n$ (for example restricted growth functions) instead of permuations of length $n$ only. We find that there is a natural way to create those posets by expressing the restricted growth functions along the maximal chains satisfying the embedding condition of posets as in the section $6$ of $\cite{MR1644460}$. From there we further note that due to satisfying the embedding conditions, those posets form a Projective Fraisse family of trees with monotone bonding epimorphisms making the corresponding inverse limit arcwise connected, dendrite and kelley as in \cite{MR1644473} and
\cite{MR1644474}. We further describe a natural way to create lattice from those trees and study the above mentioned properties investigating the confluency and monotonicity of the corresponding bonding maps. We tried to investigate the structure of the projective limits of Hasse diagrams in the scenareo of some other posets. For example, Jean-Luc Baril, Mireille Bousquet-Melou, Sergey Kirgizov, Mehdi Naima worked on Ascent lattice of dyck path in \cite{MR1644477} extending the idea of well known Stanley Order into a greedy version of that. We notice that always $D_{n}$ carrying a copy of $D_{n-1}$ (more than one copies), we can create the profinite graph $\varprojlim_{n\in N}{D_{n}}$ ($N$ being the set of all positive integers) from their Hasse diagrams and the corresponding bonding maps being confluent, the projective limit is Kelley. We get similar structure studying $p$-generalised Fribonacci sequences avoiding certain patterns in the corresponding Dyck paths as mentioned in \cite{MR1644479}. There has been a lot of study in the literature on higher dimer covers of Snake graphs of continued fractions $[a_{1},a_{2},\cdots a_{n}]$ (where $a_{1}, a_{2}\cdots a_{n}$ are positive integers) by Gregg Musiker, Nicholas Ovenhouse, Ralf Schiffler, and Sylvester W. Zhang and Amanda Burcroff as in \cite{MR1644475}\cite{MR1644476}.  We found that under a particular scenareo the collections of lattices of all $m$ dimer covers of the border strip $G[n,n]$ for $m,n \in N$ forms a Projective Fraisse family of graphs with confluent bonding epimorphism, making the projective limit Kelley. We further notice that the similar situation arrives in case of lattices of flats of Uniform matroids.
There has been a lot of study on non crossing set partitions of $[n]$, for example as in \cite{MR1644478}. We notice that in a particular case, when the number of blocks of a set partition is $k=2$, the collection of all lattices of noncrossing partitions of $[n]$, varying the $n$ upto infinity forms projective Fraisse family of graphs with monotone bonding epimorphism making the projective limit Kelley. In order to extend the idea of algebraic and combinatorial approach of Kraft, De Concini-Procesi, and Tanisaki from the Springer varieties to the regular nilpotent Hessenberg varieties, Aba Mbirika, Julianna Tymoczko worked on lattices of Hessenberg functions in \cite{MR1644490}. Regular nilpotent Hessenberg varieties have been studied extensively in \cite{MR1644480}, \cite{MR1644481}, \cite{MR1644482}, \cite{MR1644483}, \cite{MR1644484}. We notice that, the lattices of Hessenberg functions form a projective Fraisse families of graphs with confluent epomorphism, making the inverse limit Kelley. 

\section{Preliminary}

It is known from section 6 in \cite{MR1644460}
when we call a poset P, a
$\mathcal{P}$ - chain - permutational
given a subset of permutations $\mathcal{P}$
of $S_{n}$. In this work, we use the
same idea to study subset of
words that are not necessarily
permutations for example
especially when they are certain
classes of restricted growth
functions (r.g.f.) induced by pattern avoidance classes.
\begin{definition}
A poset P will be called
$\mathcal{P}$ - chain - word, where $\mathcal{P}$ is a
subset of the set of all words of
length $n$, if it is possible to label
the covering relations of P (i.e.
the edges of the Hasse diagram of
P) with numbers from 1, 2, ..., n
in such a way that along different
maximal chains of P (whose
length is necessarily $n$) the labels
form different words from P, and
every word w of P arises in this
manner. 
\end{definition}
For example, it is well
known that the Boolean lattice
$B_{n}$ is a chain word poset where
the set of words is all the words (in fact in this example, they are permutations) in one line notation in $S_{n}$.
\begin{definition}\cite {MR1544457} A restricted growth function (RGF) is a sequence $w = a_{1} ...a_{n}$ of positive integers subject to the
	restrictions 
	\begin{enumerate}
		\item $a_{1} =1$. 
		\item For $i\geq 2$,  $a_{i} \leq 1+max\{a_{1},...,a_{i-1}\}$
	\end{enumerate}

 \end{definition}

In \cite {MR1544457} a partition of $[n]$ is written as $\sigma = B_{1}/\cdots /B_{k}\vdash S$ where the subsets $B_{i}$ are called blocks. We use the notation $\Pi_{n} =\{\sigma : \sigma \vdash [n]\}$.
In order to connect set partitions restricted growth functions as in \cite {MR1544457} we need the elements of $\Pi_{n}$ in standard form. To define pattern avoidance in this setting, suppose $\sigma = B_{1}/\cdots /B_{k} \in \Pi_{n}$ and $S \subseteq [n]$. Then $\sigma$ has a corresponding sub partition $\sigma^{'}$ whose blocks are the nonempty intersections $B_{i} \bigcap S$. For example, if $\sigma$ $=16/23478/5 \vdash [8]$ and $S =\{2, 4, 6, 7\}$ then $\sigma^{\prime} = 247/6$. We standardize a set partition by replacing the smallest element by $1$, the next smallest by $2$, and so on. So the standardization of $\sigma^{\prime}$ above is $124/3$. Given two set partitions $\sigma$ and $\pi$, we say that $\sigma$ contains $\pi$ as a pattern if there is a sub partition of $\sigma$ which standardizes to $\pi$. Otherwise, we say that $\sigma$ avoids $\pi$. Continuing our example, we have already shown that $\sigma = 16/23478/5$ contains $124/3$. But $\sigma$ avoids $12345/6$ because the only block of $\sigma$ containing five elements also contains the largest element in $\sigma$, so there can be no larger element in a separate block. As in  we let $\Pi_{n}(\pi) =\{\sigma \in \Pi_{n} : \sigma$ avoids $\pi\}$.
\begin{definition}
	We say $\sigma = B_{1}/\cdots /B_{k} \in \Pi_{n}$ is in standard form if min $B_{1} < \cdots <$ min $B_{k}$. Thus it follows that min $B_{1} =1$.
\end{definition}

We assume all partitions in $\Pi_{n}$ are written in the standard form. Associate with $\sigma\in\Pi_{n}$ the word $w(\sigma) = a_{1}\cdots a_{n}$ where $a_{i} =j$ if and only if $i\in B_{j}$. For example  $w(16/23478/5) =12223122.$
Let, $\Pi_{n,k}$ be the set of all words in $\Pi_{n}$ with exactly $k$ many blocks. $R_{n} =\{w : w$ is an RGF of length $n\}$. Let, $R_{n,k} =\{w : w$ is an RGF of length $n$ with maximal letter $k\}$.

\begin{theorem} The characterizations below follow from theorem $1.2$ in  \cite{MR1544457}  
and \cite{MR1644460} 
\begin{enumerate}
\item[i.] $R_{n,k}(1/2/3)= 1$, if $k=1$ and this is equal to $\{w\in R_{n,k}:w$ consists of only 1 s and 2 s$\}$, if $(k=2)$. $R_{n,k}(1/2/3)= \Phi$, if $k>2$.
\item[ii.] $R_{n,k}(1/23) =\{w\in R_n : w $ is obtained  by inserting  a  single $1$ into  a  word of the form  $1^{l}23\cdots  k$ for some $l\geq 0\}$
\item[iii.] $R_{n,k}(13/2)$  = $\{w\in R_{n}: w$ is layered ie $w = 1^{l_{1}}2^{l_{2}}\cdots  k^{l_{k}}\}$, for some positive integers $l_{1}, l_{2},\cdots  l_{k}$.
\item[iv.]  $R_{n,k}(12/3) =\{w \in R_n : w$ has initial run $1...k$ and $a_{k+1}=\cdots  a_{n} \leq k\}$.
\item[v.] $R_{n,k}(123) =\{w\in R_{n,k}:w$ has no element repeated more than twice$\}$
\end{enumerate}

\end{theorem}
We state below the theorem $1.2$ from \cite{MR1544457}
\begin{theorem}
\begin{enumerate}
\item[i.] $R_{n}(1/2/3)= \{w\in R_{n}:w$ consists of only 1 s and 2 s only $\}$. 
\item[ii.] $R_{n}(1/23) =\{w\in R_n : w $ is obtained  by inserting  a  single $1$ into  a  word of the form  $1^{l}23\cdots  m$ for some $l\geq 0 and m\geq 1\}$
\item[iii.] $R_{n}(13/2)$  = $\{w\in R_{n}: w$ is layered ie $w = 1^{l_{1}}2^{l_{2}}\cdots  m^{n_{m}}\}$, for some positive integers $l_{1}, l_{2},\cdots  l_{m}$.
\item[iv.] $R_{n}(12/3) =\{w \in R_n : w$ has initial run $1...m$ and $a_{m+1}=\cdots  a_{n} \leq m\}$.
\item[v.] $R_{n}(123) =\{w\in R_{n}:w$ has no element repeated more than twice $\}$
\end{enumerate}

Now we provide some backgrounds from \cite{MR1644473} and \cite{MR1644474} on Projective Fraisse family of trees and graphs here.
 \begin{definition}A graph is an ordered pair $A = (V(A),E(A))$, where $E(A) \subseteq V(A)^{2}$ is a reflexive
 and symmetric relation on $V(A)$. The elements of $V(A)$ are called vertices of graph
 $A$ and elements of $E(A)$ are called edges. Given vertices $a$ and $b$ in a graph $A$, where
 $a\neq b$, we will use $\langle a,b \rangle$ to denote an edge in $A$. 
\end{definition}
Note that since $E(A)$ is reflexive there
 is always an edge from a vertex in $A$ to itself, such an edge will be called a degenerate
 edge. When it is clear from context that $x$ is a vertex, we will use $x \in A$ for $x \in V(A)$.
 \begin{definition}Given two graphs $A$ and $B$ a function $f:V(A) \mapsto V(B)$ is a homomorphism if it
 maps edges to edges, i.e. $\langle a,b \rangle \in E(A)$ implies $\langle f(a),f(b)\rangle \in E(B)$. 
\end{definition}
A homomorphism
 $f$ is an epimorphism between graphs if it is moreover surjective on both vertices and
 edges. An isomorphism is an injective epimorphism.
 
 By a rooted graph we mean a graph $G$ with a distinguished vertex $r_{G}$. An epimorphism $f: G \mapsto H$, between rooted graphs, is an epimorphism between corresponding
 (unrooted) graphs, such that $f(r_{G}) = r_{H}$
\begin{definition} 
 	A topological graph (rooted topological graph) $K$ is a graph $(V(K),E(K))$
 (respectively, rooted graph $(V(K),E(K),r_{K})$, where $r_{K}$ is the distinguished vertex),
 whose domain $V(K)$ is a $0$-dimensional, compact, second-countable (thus metrizable)
 space and $E(K)$ is a closed, reflexive and symmetric subset of $V (K)^{2}$. 
\end{definition}
We require that
 all epimorphisms between topological graphs are continuous.

\begin{definition}
	 A tree $T$ is a finite graph such that for every two distinct vertices
$a, b \in T$ there is a unique finite sequence $v_{0} = a,v_{1},\cdots,v_{n} = b$ of vertices in $T$ such
that for every $i \in\{ 0\}\bigcup [n] $, we have $\langle v_{i},v_{i+1} \rangle \in E(T)$ and $v_{j} \neq v_{i}$, for $i\neq j$.
\end{definition}
\begin{definition}
 	 Let $F$ be a nonempty family of finite graphs or a family of finite rooted
 graphs with a fixed family of epimorphisms among the structures in $F$. We say that $F$
 is a projective Fraisse family if
 \begin{enumerate}
 	\item $F$ is countable up to isomorphism, that is, any sub-collection of pairwise non
 	isomorphic structures of $F$ is countable
 	
 	\item  Epimorphisms in $F$ are closed under composition and each identity map is in $F$
 	
 	\item For $B,C \in F,  \exists D \in F$ and epimorphisms $f$ and $g$ in $F$ such that
 	$f : D \mapsto B$ and $g: D\mapsto C$ and
 	
 	\item  for every two epimorphisms $f : B \mapsto A$ and $g: C \mapsto A$ in $F$ there exist epimor
 	phisms $f_{0}: D \mapsto B$ and $g_{0}: D \mapsto C$ in $F$ such that $f \circ f_{0} = g\circ g_{0}$.
 \end{enumerate}
 \end{definition}
 	
 	Notation: Given a sequences $F_{n}$ of finite (rooted) graphs and epimorphisms
 	$\alpha_{n}: F_{n+1} \mapsto F_{n}$, which are called bonding maps, we denote the inverse sequence by
 	({$F_{n},\alpha_{n}$}). For $m > n$, we let given
 	an inverse sequence ({$F_{n},\alpha_{n}$}), the associated inverse limit space is the subspace of the
 	product space $\Pi{F_{i}}$ determined by $\{(x_1,x_2,\cdots x_i) \in F_i\}$ and $x_i = \alpha_{i}(x_{i+1})\}$ and is
 	denoted as $F = \varprojlim {F_{n}}
 	$. Further, we denote the canonical projection from the
 	inverse limit space $F$ to the $n-th$ factor space $F_{n}$ by $\alpha^{\infty}_{n}$. Note that if $x,y \in F $,
 	then the metric $d(x,y) = \sum_
 	{i=1}^{\infty}
 	\frac{d_{i}(xi,yi)}{2^{i}}$, where $d_{i}$ is the discrete metric on $F_{i}$, is
 	compatible with the topology of the product space.\\
 	
 	The family $F$ of finite graphs or finite rooted graphs and epimorphisms is enlarged
 	to a family $F^{\omega}$ which includes all (rooted) topological graphs obtained as inverse limits
 	of (rooted) graphs in $F$ with bonding maps from the family of epimorphisms. If $G =
 	\varprojlim
 	{G_{n},\alpha_{n}} \in F^{\omega}$ and $a = (a_{n})$ and $b = (b_{n})$ are elements of $G$ then $\langle a,b\rangle$ is an edge
 	in $G$ if and only if for each $n$, $\langle a_{n},b_{n}\rangle$ is an edge in $G_{n}$.

\begin{definition} Given two topological graphs G and H, an epimorphism
$f : G \mapsto H$ is called confluent if for every closed connected subset
$Q$ of $H$ and every component $C$ of $f^{-1}(Q)$ we have $f(C) = Q$. Equivalently,
if for every closed connected subset $Q$ of $H$ and every vertex
$a \in G$ such that $f(a) \in Q$, there is a connected set $C$ of $G$ such that
$a \in C$ and $f(C) = Q$. It is called monotone if preimages of closed
connected sets are connected.
\end{definition}
 Thus every monotone epimorphism is
confluent.
 If $f : X \mapsto Y$ and $g : Y \mapsto Z$ are confluent epimorphisms
between topological graphs, then $g \circ f : X \mapsto Z$ is a confluent
epimorphism.
  
\end{theorem}
		\section{Chain word posets and projective Fraisse families of trees from set partitions and pattern avoidence classes:}

We first notice that given any positive integer $n$, there is a natural way to create a rooted $n$-ary tree structure (viewed as a poset) labeling the top edge from the root always as $1$ (since the first digit of a restricted growth function is always $1$) and along different maximal chains of the tree (whose
length is necessarily $n$) the labels form different words from $R_{n}$, and every word $w$ of $R_{n}$ arises in this
manner, displaying any $R_{n}$ as a $\mathcal{P}$-chain poset, where $\mathcal{P}=R_{n}$. Also, we see a natural embedding of $R_{n}$ inside any $R_{m}$, where $m>n$, since any R.G.F. of length $n$ can be extended to that of length $n+1$ by inserting the digit $1$ at the end. And that is why backwardly, we can find a natural projection of the Hasse diagram of that poset $R_{m}$ upon $R_{n}$ as an epimorphism of graphs where the top part of the rooted tree $R_{m}$ containing the  copy of $R_{n}$ is identically mapped onto $R_{n}$ and the rest of all the edges below any vertex of that copy of $R_{n}$ can be mapped onto the loop around that specific vertex in $R_{n}$, creating a projective Frasse Family of trees (and hence a poset, if we consider the elements of the posets as the vertices of the respective trees) with  monotone epimorphisms.
Look at the figure in case of type B set partitions utilizing generalised or signed restricted growth functions.

 We note that for any set pattern $v$ (and/or multiple patterns $v$,$w$, etc.), the similar as above happens in case of all set partitions of $[n]$ avoiding $v$ considering the corresponding restricted growth functions $R_{n}(v)$,(and/or multiple patterns $R_{n}(v,w)$, etc.) varying $n$ over the set of all positive integers. 

If we call the corresponding projective Fraisse limit of those trees as \textbf{T}, the by theorem 2.18 in \cite{MR1644473}, T is heriditarily unicoherent as a continum (as in the definition of the metric topology in the above notation in preliminary section), and arcwise connected (by  Proposition 3.10 \cite{MR1644473}) and is a dentrite by definition 2.10 and corollary 3.11 in \cite{MR1644473}. By definition 5.10 and proposition 5.12 in \cite{MR1644473}, \textbf{T} also smooth dendroid. In addition, by proposition 6.5 and definition 6.3 in \cite{MR1644473} is also Kelley. Also by section 6 in (\cite{MR1644473}), they have lifting property.\\

We further notice that, we can extend those trees as a lattice in a very natural way tagging an $\epsilon$ in the bottom of those trees and connecting that $\epsilon$ via an edge with each vertex in the bottom level of that tree. We denote the corresponding lattices from the pattern avoidance classes as $R_{n,L}(\sigma)$ where $\sigma$ is a set of patterns. And we note that those monotone bonding epimorphisms can be extended naturally from $R_{n+1, L}$ onto $R_{n,L}$ by mapping the copy of corresponding tree from $R_{n,L}$ in the top of $R_{n+1,L}$ identically, and mapping all the edges from the bottom vertex of the copy of $R_{n,L}$ inside the domain to the corresponding loops and the vertex $\epsilon$ onto $\epsilon$.
Also, it may be of further interest to study these Projective Fraisse Family of trees (and lattices) varying both $n$ and $k$ and providing a fixed $k$ varying the $n$ only. We discuss specific examples in all these scenario. Often times, we consider $n>k$, instead of $n\geq k$, to avoid triviality for $n=k$.
\begin{example}
\item Consider the set of all $R_{n,k}(1/23)$. \\
Set up $A=\{(n,k): n,k \in N, n> k\}$. For a given fixed $n$, if we increase $k$ keeping $1<k<n$, then we notice that $R_{n,k+1}(1/23)$ contains a copy of $R_{n,k}(1/23)$, where $1<k<n$, we can define an epimorphism of graphs $\phi_{n,k+1}:R_{n,k+1}(1/23)\rightarrow R_{n,k}(1/23)$ as follows, the copy of $R_{n,k}(1/23)$ sitting inside $R_{n,k+1}(1/23)$ in the left is identically mapped onto that. The right most maximal chain $1^{n-(k+1)}23\cdots(k+1)$ inside $R_{n,k+1}(1/23)$ falls upon $1^{n-(k)}23\cdots k$ inside $R_{n,k}(1/23)$ creating a Projective Fraisse family of trees and confluent bonding epimorphisms (although, note that this bonding maps are not necessarily monotone, since the pre-image of the right bottom vertex in the range graph is not necessarily connected), while we define $(n_{1}, k_{1})\geq (n_{2}, k_{2})$, with $n_{1}\geq n_{2}, k_{1}\geq k_{2}$. The family is really a Projective Fraisse family, since because if we consider an $R_{n,k}(1/23)$ and $R_{n+1,k}(1/23)$, then $R_{n+1,k}(1/23)$ contains copy of $R_{n,k}(1/23)$ in it's bottom. In this case, if we only vary the $n$ upto infinity fixing $k$, we notice that the bonding maps are monotone. Extending the bonding maps naturally in the corresponding lattice structures varying both $n$ and $k$, we get a projetcive Fraisse family of graphs (from the Hasse diagrams of those lattices) with confluent epimorphisms. And in the lattice case, if we fix $k$ and make the $n$ vary upto infinity, the corresponding bonding maps are again monotone. By proposition 6.5 and definition 6.3 in $\cite{MR1644473}$, the corresponding projective limit, while both $n, k$ vary, is Kelley (in case of tree and Lattice both). When, we fix $k$, make the $n$ going to infinity by proposition 3.10 in \cite{MR1644473}, the corresponding inverse limit is arcwise connected (in case of tree and lattice both)
\item Consider the set of all $R_{n,k}(12/3)$\\
 We notice first that these graphs has a very specific shape due to the initial run in the begining which is as follows: always we have only one path of length $k$ from the root and from the tip of the end of that path, there are $k$ many paths of length $n-k$. Precisely that's what is the shape of the tree. We know that in order to have a valid $(n,k)$ in this set up we must need $n>k$. So, consider the set $A=\{(n,k): n,k \in \mathcal {N}, n> k\}$, note that if we define $(n_{1}, k_{1})\geq (n_{2}, k_{2})$, with $n_{1}\geq n_{2}, k_{1}\geq k_{2}$ and additionally this time $(n_{1}-n_{2})\geq (k_{1}-k_{2})$. In that case, we can define epimorphism of graphs $\phi:R_{(n_{1}, k_{1})}(12/3)\mapsto R_{(n_{2}, k_{2})}(12/3)$ as follows: The $k_{1}-k_{2}$ many edges in the top of the domain graph are mapped onto the loop at the root in the destination graph. Next we see, the rest of the domain graph contains a copy of the destination, so we map that portion identically on the range graph and at the bottom of each branch the $(n_{2}-k_{2})-(n_{1}-k_{1})$ many edges are mapped onto the corresponding loop at the corresponding end vertex in the domain. Then the rest of the $k_{2}-k_{1}$ many branches are mapped on to the left most branch, wheras $(n_{2}-k_{2})-(n_{1}-k_{1})$ many edges in the bottom of each such branch are mapped on to the loop at the corresponding end vertex.  All these creates a Projective Fraisse Family of trees with confluent epimorphisms (considering the above preorder) that are not necessarily monotone. In this case, if we only vary the $n$ upto infinity fixing $k$, we notice that the bonding maps are monotone. As in the previous example we can extend the confluency (when n,k both varies) and monotonicity (when, n vary for a fixd $k$) of the bonding epimorphisms. By proposition 6.5 and definition 6.3 in $\cite{MR1644473}$, the corresponding projective limit is while both $n, k$ vary, is Kelley (in case of tree and Lattice both). When, we fix $k$, make the $n$ going to infinity by proposition 3.10 in \cite{MR1644473}, the corresponding inverse limit is arcwise connected (in case of tree and lattice both)

\item One can extend this idea for standard type B set partitions as in \cite{MR1644458}  
,\cite{MR1644469} ,\cite{MR1644470},\cite{MR1644471} as well.  Note that, if we denote these by $R^{B}_{n,k,L}$, then we can label the maximal chains of those lattices also by using the signed restricted growth functions (and the corresponding bijections of standard type B partitions with them) as in \cite{MR1644470}, \cite{MR1644471}. We notice that the Hasse diagram of $R^{B}_{n,k,L}$ contains a copy of that of $R^{B}_{(n-1),k,L}$ to it's left at the at the top, creating again the projection in a natural way as in the above examples, giving a projective Fraisse family of graphs (induced by the trees in a natural way as before in case of set partitions of $[n]$) with monotone epimorphisms. If we ignore the lattice structure, omitting the '$\epsilon$' and the edges connecting to $\epsilon$, then we see again the corresponding projective limit is arcwise connected, Kelley and a smooth dendroid, satisfying lifting property by \cite{MR1644473}. Some figures are given below for type B.
\newpage
\adjustbox{scale=0.45,center}
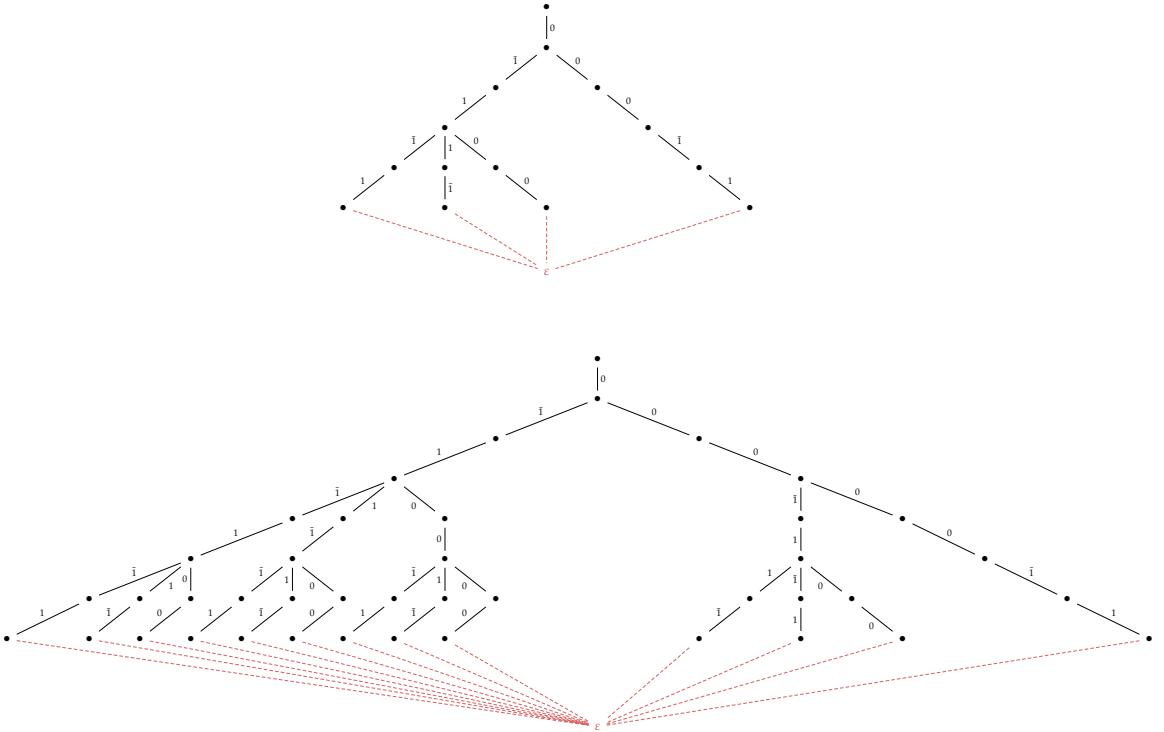
{%
\begin{tikzcd}
	&&&&&&&&&&& \bullet \\
	&&&&&&&&&&& \bullet \\
	&&&&&&&&&& \bullet && \bullet \\
	&&&&&&&&& \bullet &&&& \bullet \\
	&&&&&&&& \bullet & \bullet & \bullet &&&& \bullet \\
	&&&&&&& \bullet && \bullet && \bullet &&&& \bullet \\
	\\
	&&&&&&&&&&& \textcolor{rgb,255:red,214;green,92;blue,92}{{\Huge\varepsilon}} \\
	\\
	\\
	&&&&&&&&&&&& \bullet \\
	&&&&&&&&&&&& \bullet \\
	&&&&&&&&&& \bullet &&&& \bullet \\
	&&&&&&&& \bullet &&&&&&&& \bullet \\
	&&&&&& \bullet & \bullet && \bullet &&&&&&& \bullet && \bullet \\
	&&&& \bullet && \bullet &&& \bullet &&&&&&& \bullet &&&& \bullet \\
	&& \bullet & \bullet & \bullet & \bullet & \bullet & \bullet & \bullet & \bullet & \bullet &&&&& \bullet & \bullet & \bullet &&&&& \bullet \\
	\bullet && \bullet & \bullet & \bullet & \bullet & \bullet & \bullet & \bullet & \bullet &&&&& \bullet && \bullet && \bullet &&&&&& \bullet \\
	\\
	\\
	&&&&&&&&&&&& \textcolor{rgb,255:red,214;green,92;blue,92}{{\Huge\varepsilon}}
	\arrow["0", no head, from=1-12, to=2-12]
	\arrow["{\bar{1}}"', no head, from=2-12, to=3-11]
	\arrow["0", no head, from=2-12, to=3-13]
	\arrow["1"', no head, from=3-11, to=4-10]
	\arrow["0", no head, from=3-13, to=4-14]
	\arrow["{\bar{1}}"', no head, from=4-10, to=5-9]
	\arrow["1", no head, from=4-10, to=5-10]
	\arrow["0", no head, from=4-10, to=5-11]
	\arrow["{\bar{1}}", no head, from=4-14, to=5-15]
	\arrow["1"', no head, from=5-9, to=6-8]
	\arrow["{\bar{1}}", no head, from=5-10, to=6-10]
	\arrow["0", no head, from=5-11, to=6-12]
	\arrow["1", no head, from=5-15, to=6-16]
	\arrow[draw={rgb,255:red,214;green,92;blue,92}, dashed, no head, from=6-8, to=8-12]
	\arrow[draw={rgb,255:red,214;green,92;blue,92}, dashed, no head, from=6-10, to=8-12]
	\arrow[draw={rgb,255:red,214;green,92;blue,92}, dashed, no head, from=6-12, to=8-12]
	\arrow[draw={rgb,255:red,214;green,92;blue,92}, dashed, no head, from=6-16, to=8-12]
	\arrow["0", no head, from=11-13, to=12-13]
	\arrow["{\bar{1}}"', no head, from=12-13, to=13-11]
	\arrow["0", no head, from=12-13, to=13-15]
	\arrow["1"', no head, from=13-11, to=14-9]
	\arrow["0", no head, from=13-15, to=14-17]
	\arrow["{\bar{1}}"', no head, from=14-9, to=15-7]
	\arrow["1", no head, from=14-9, to=15-8]
	\arrow["0"', no head, from=14-9, to=15-10]
	\arrow["{\bar{1}}"', no head, from=14-17, to=15-17]
	\arrow["0", no head, from=14-17, to=15-19]
	\arrow["1"', no head, from=15-7, to=16-5]
	\arrow["{\bar{1}}"', no head, from=15-8, to=16-7]
	\arrow["0"', no head, from=15-10, to=16-10]
	\arrow["1"', no head, from=15-17, to=16-17]
	\arrow["0", no head, from=15-19, to=16-21]
	\arrow["{\bar{1}}"', no head, from=16-5, to=17-3]
	\arrow["1", no head, from=16-5, to=17-4]
	\arrow["0"', no head, from=16-5, to=17-5]
	\arrow["{\bar{1}}"', no head, from=16-7, to=17-6]
	\arrow["1"', no head, from=16-7, to=17-7]
	\arrow["0"', no head, from=16-7, to=17-8]
	\arrow["{\bar{1}}"', no head, from=16-10, to=17-9]
	\arrow["1"', no head, from=16-10, to=17-10]
	\arrow["0"', no head, from=16-10, to=17-11]
	\arrow["1"', no head, from=16-17, to=17-16]
	\arrow["{\bar{1}}"', no head, from=16-17, to=17-17]
	\arrow["0"', no head, from=16-17, to=17-18]
	\arrow["{\bar{1}}", no head, from=16-21, to=17-23]
	\arrow["1"', no head, from=17-3, to=18-1]
	\arrow["{\bar{1}}"', no head, from=17-4, to=18-3]
	\arrow["0"', no head, from=17-5, to=18-4]
	\arrow["1"', no head, from=17-6, to=18-5]
	\arrow["{\bar{1}}"', no head, from=17-7, to=18-6]
	\arrow["0"', no head, from=17-8, to=18-7]
	\arrow["1"', no head, from=17-9, to=18-8]
	\arrow["{\bar{1}}"', no head, from=17-10, to=18-9]
	\arrow["0"', no head, from=17-11, to=18-10]
	\arrow["{\bar{1}}"', no head, from=17-16, to=18-15]
	\arrow["1"', no head, from=17-17, to=18-17]
	\arrow["0"', no head, from=17-18, to=18-19]
	\arrow["1", no head, from=17-23, to=18-25]
	\arrow[draw={rgb,255:red,214;green,92;blue,92}, dashed, no head, from=18-1, to=21-13]
	\arrow[draw={rgb,255:red,214;green,92;blue,92}, dashed, no head, from=18-3, to=21-13]
	\arrow[draw={rgb,255:red,214;green,92;blue,92}, dashed, no head, from=18-4, to=21-13]
	\arrow[draw={rgb,255:red,214;green,92;blue,92}, dashed, no head, from=18-5, to=21-13]
	\arrow[draw={rgb,255:red,214;green,92;blue,92}, dashed, no head, from=18-6, to=21-13]
	\arrow[draw={rgb,255:red,214;green,92;blue,92}, dashed, no head, from=18-7, to=21-13]
	\arrow[draw={rgb,255:red,214;green,92;blue,92}, dashed, no head, from=18-8, to=21-13]
	\arrow[draw={rgb,255:red,214;green,92;blue,92}, dashed, no head, from=18-9, to=21-13]
	\arrow[draw={rgb,255:red,214;green,92;blue,92}, dashed, no head, from=18-10, to=21-13]
	\arrow[draw={rgb,255:red,214;green,92;blue,92}, dashed, no head, from=18-15, to=21-13]
	\arrow[draw={rgb,255:red,214;green,92;blue,92}, dashed, no head, from=18-17, to=21-13]
	\arrow[draw={rgb,255:red,214;green,92;blue,92}, dashed, no head, from=18-19, to=21-13]
	\arrow[draw={rgb,255:red,214;green,92;blue,92}, dashed, no head, from=18-25, to=21-13]
\end{tikzcd}
}

\end{example}
Note: We further notice that $R_{n}(1/2/3), R_{n}(13/2), R_{n}(123), R_{n}(1/23)$ being binary trees, the corresponding lattices are distributive with or without $k$ and sublattice and products of distributive lattices, being so, the corresponding projective limits are infinite distributive latices. On the other hand, $R_{n,L}(12/3)$ and $R^{B}_{n,k,L}$ are not distributive, where as induced by set partitions over $[n]$ and type B partitions over $\langle n\rangle$, they are ranked. By section 6 in \cite{MR1644473}, all the families of graphs in the above example have lifting property. The way the lattices are constructed tagging an "$\epsilon$" only at the end of the corresponding tree, it follows from lemma 6.2 in \cite{MR1644473}, the topological graph $G$ induced by taking their projective limit is arcwise connected. Also, the bonding maps in between the corresponding trees being monotone, when we extend those rooted tree structures in to a lattice, by tagging the "$\epsilon$" in a natural way, the monotonicity in the bonding maps preserve (by lemma 3.4 in \cite{MR1644473} and that is why by proposition 3.9, the corresponding projections from the inverse limit $G$ on to each finite factor 
stays as monotone. We notice that from the outer bounds (and many other alternative ways) of those lattices, we can extract systems of cycles and redefine the bonding maps in between those cycles to construct graph solenoid as in \cite{MR1644474}. For example: consider the family of cycles from the outer boundries of $R_{n,L}(1/23)$ for $n=2, 2^{2}, 2^{3},\cdots$, ignoring the top edge from the root. It is not hard to see that wrapping the circle of the outer boundary of $R_{4,L}(1/23)$ along that of $R_{2,L}(1/23)$ ignoring the top edge from the root (mapping two consequtive vertices in the domain onto the same vertex, then next two consequtive vertices on the next same vertices and so on in the range) we get a confluent epimorphism where the winding number is $2$. Continuing this process for all $n=2^{n}$, we get an inverse system of cycles and confluent epimorphisms $(C_{n}, p_{n})$, making the inverse limit of graph sole-noid as in \cite{MR1644474}.
$R_{2^{n}L}(1/23)$, we get an inverse sequence of cycles $(C_{n}, p_{n})$, where each $p_{n}$ is a confluent epimorphism.

	 \section{projective Frasse families of graphs from other lattices}
\begin{enumerate}
\item Here we provide some information from \cite{MR1644477} first.
Let us first recall that a Dyck path $P$ is a 2-dimensional path starting at the
origin $(0,0)$, consisting of up steps $U=(1,1)$ and down steps $D=(1,-1)$, that ends on the
$x$-axis and never goes strictly below the $x$-axis. The size of P is the number n of up steps. We
denote by $D_{n}$, the set of such paths. An ascent of a path $P$ is a maximal, non-empty sequence
of consecutive up steps. A descent is defined similarly using down steps. A factor of $P$ is a
non-empty sequence of consecutive steps. A peak is a factor $UD$, while a valley is a factor
$DU$. There exists on $D_{n}$ a classical order, called the Stanley order , for which $P$ is less than or equal to $Q$ if it lies (weakly) below $Q$. By this,
we mean that for any $l$, the prefix of $P$ of length $l$ contains at most as many up steps as the
prefix of length $l$ of $Q$. A path $Q$ covers a path $P$ in this order if and only if $Q$ is obtained by
replacing a valley $DU$ of $P$ into a peak $UD$. In \cite{MR1644477} they consider a greedy version of this order, which they call the ascent order. It
is described by its cover relations: they say that $Q$ covers $P$ (denoted $P <·Q$) if there exists in $P$
a factor $DU^{k}D$ such that Q is obtained from P by replacing this factor by $U^{k}DD$. Observe
that in this case $P$ lies below $Q$. This relation being irreflexive and anti-symmetric,
its transitive closure is thus an order relation on $D_n$, denoted $\leq$, and one checks that
the cover relations of this order are indeed those described above.  The
corresponding poset is denoted by $D_{n}$.

We can define an epimorphism of graphs and posets $f_{n}:D_{n}\mapsto D_{n-1}$ as follows $f_{n}(P)$ is obtained by deleting the first peak $UD$ of $P$, by utilising the fact that $D_{n}$ contains a copy of $D_{n-1}$ to it's bottom right (there are other copies of $D_{n-1}$ inside $D_{n}$). It follows that each $f_{n}$ is monotone, as each vertex fiber is connected. Hence, by definition 6.3 and proposition 6.5 in \cite{MR1644473}, the corresponding inverse limit $D$ is Kelley. Each $D_{n}$ is a non distributive lattice. Hence, the corresponding inverse limit $D$ is non distributive. The necessary figure is given below.

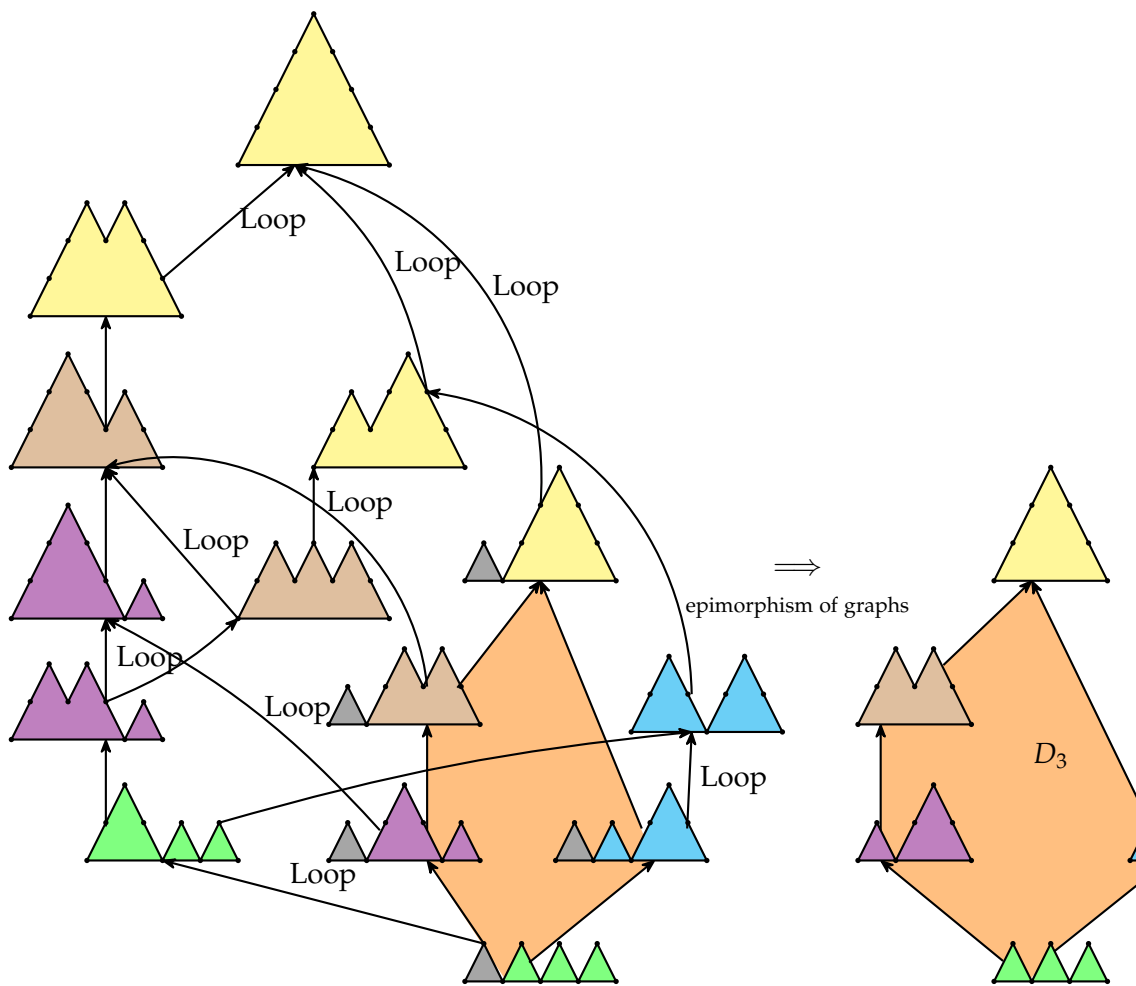
\begin{tikzpicture}[
  thick,
  >={Stealth[round]},
  dot/.style={circle,fill=black,inner sep=0pt,minimum size=2.0pt},
  cYellow/.style={fill=yellow!50},
  cOrange/.style={fill=orange!50},
  cBlue/.style={fill=cyan!50},
  cBrown/.style={fill=brown!50},
  cPurple/.style={fill=violet!50},
  cGreen/.style={fill=green!50},
  cBlack/.style={fill=gray!70},
  cYellowEdge/.style={draw=yellow!60},
  cOrangeEdge/.style={draw=orange!70},
  cBlueEdge/.style  ={draw=cyan!70},
  cBrownEdge/.style ={draw=brown!70},
  cPurpleEdge/.style={draw=violet!70},
  cGreenEdge/.style ={draw=green!70},
  cBlackEdge/.style ={draw=gray!70},
]

% ------------------------------------------------
% Helper: draw one Dyck path polygon with dots
% #1 = color style, #2 = path " (..,..) -- ... ", #3 = list "(..,..),..."
% ------------------------------------------------
\newcommand{\DyckPic}[3]{%
  \filldraw[#1] #2 -- cycle;
  \foreach \p in {#3} {
    \node[dot] at \p {};
  }
  % anchor in the middle of the picture
  \path (current bounding box.center) coordinate (-center);
}

% ------------------------------------------------
% D3 pics  (length 6)
% vertices were computed from the Dyck words with step (dx,dy)=(0.25,0.5)
% ------------------------------------------------

% UDUDUD  (bottom, sawtooth with 3 peaks)
\tikzset{
  Dthree-UDUDUD/.pic={
    \DyckPic{cGreen}%
      {(0,0) -- (0.25,0.5) -- (0.5,0) -- (0.75,0.5) -- (1,0) -- (1.25,0.5) -- (1.5,0)}%
      {(0,0),(0.25,0.5),(0.5,0),(0.75,0.5),(1,0),(1.25,0.5),(1.5,0)}
  }
}

% UDUUDD  (lower left)
\tikzset{
  Dthree-UDUUDD/.pic={
    \DyckPic{cPurple}%
      {(0,0) -- (0.25,0.5) -- (0.5,0) -- (0.75,0.5) -- (1,1) -- (1.25,0.5) -- (1.5,0)}%
      {(0,0),(0.25,0.5),(0.5,0),(0.75,0.5),(1,1),(1.25,0.5),(1.5,0)}
  }
}

% UUDUDD  (middle left)
\tikzset{
  Dthree-UUDUDD/.pic={
    \DyckPic{cBrown}%
      {(0,0) -- (0.25,0.5) -- (0.5,1) -- (0.75,0.5) -- (1,1) -- (1.25,0.5) -- (1.5,0)}%
      {(0,0),(0.25,0.5),(0.5,1),(0.75,0.5),(1,1),(1.25,0.5),(1.5,0)}
  }
}

% UUDDUD  (middle right)
\tikzset{
  Dthree-UUDDUD/.pic={
    \DyckPic{cBlue}%
      {(0,0) -- (0.25,0.5) -- (0.5,1) -- (0.75,0.5) -- (1,0) -- (1.25,0.5) -- (1.5,0)}%
      {(0,0),(0.25,0.5),(0.5,1),(0.75,0.5),(1,0),(1.25,0.5),(1.5,0)}
  }
}

% UUUDDD  (top pyramid)
\tikzset{
  Dthree-UUUDDD/.pic={
    \DyckPic{cYellow}%
      {(0,0) -- (0.25,0.5) -- (0.5,1) -- (0.75,1.5) -- (1,1) -- (1.25,0.5) -- (1.5,0)}%
      {(0,0),(0.25,0.5),(0.5,1),(0.75,1.5),(1,1),(1.25,0.5),(1.5,0)}
  }
}
%----------------------------------------
%%%%%%%%%%%%%%%%%%%%%%%%%%%%%%%%%%%%%%%%%
%            D_4 Vertices
%%%%%%%%%%%%%%%%%%%%%%%%%%%%%%%%%%%%%%%%%
%----------------------------------------

\tikzset{
  D4-UDUDUDUD/.pic={
    % Black triangle for the UD prefix
    \filldraw[cBlack] (0,0) -- (0.25,0.5) -- (0.5,0) -- cycle;
    \node[dot] at (0,0) {};
    \node[dot] at (0.25,0.5) {};
    \node[dot] at (0.5,0) {};
    % Green continuation (UDUDUD shifted)
    \filldraw[cGreen] (0.5,0) -- (0.75,0.5) -- (1,0) -- (1.25,0.5) -- (1.5,0) -- (1.75,0.5) -- (2,0) -- cycle;
    \foreach \p in {(0.75,0.5),(1,0),(1.25,0.5),(1.5,0),(1.75,0.5),(2,0)} {
      \node[dot] at \p {};
    }
    \path (current bounding box.center) coordinate (-center);
  }
}

% 6. UUDDUDUD - starts with UU (intrinsic green)
% Heights: 1,2,1,0,1,0,1,0
\tikzset{
  D4-UUDDUDUD/.pic={
    \DyckPic{cGreen}%
      {(0,0) -- (0.25,0.5) -- (0.5,1) -- (0.75,0.5) -- (1,0) -- (1.25,0.5) -- (1.5,0) -- (1.75,0.5) -- (2,0)}%
      {(0,0),(0.25,0.5),(0.5,1),(0.75,0.5),(1,0),(1.25,0.5),(1.5,0),(1.75,0.5),(2,0)}
  }
}

% 8. UUDUDDUD - starts with UU (intrinsic green)
% Heights: 1,2,1,2,1,0,1,0
\tikzset{
  D4-UUDUDDUD/.pic={
    \DyckPic{cPurple}%
      {(0,0) -- (0.25,0.5) -- (0.5,1) -- (0.75,0.5) -- (1,1) -- (1.25,0.5) -- (1.5,0) -- (1.75,0.5) -- (2,0)}%
      {(0,0),(0.25,0.5),(0.5,1),(0.75,0.5),(1,1),(1.25,0.5),(1.5,0),(1.75,0.5),(2,0)}
  }
}

% 9. UUDUDUDD - starts with UU (intrinsic green)
% Heights: 1,2,1,2,1,2,1,0
\tikzset{
  D4-UUDUDUDD/.pic={
    \DyckPic{cBrown}%
      {(0,0) -- (0.25,0.5) -- (0.5,1) -- (0.75,0.5) -- (1,1) -- (1.25,0.5) -- (1.5,1) -- (1.75,0.5) -- (2,0)}%
      {(0,0),(0.25,0.5),(0.5,1),(0.75,0.5),(1,1),(1.25,0.5),(1.5,1),(1.75,0.5),(2,0)}
  }
}

% --- PURPLE FIBER (maps to UDUUDD) ---

% 2. UDUDUUDD - UD + UDUUDD
\tikzset{
  D4-UDUDUUDD/.pic={
    % Black triangle for UD prefix
    \filldraw[cBlack] (0,0) -- (0.25,0.5) -- (0.5,0) -- cycle;
    \node[dot] at (0,0) {};
    \node[dot] at (0.25,0.5) {};
    \node[dot] at (0.5,0) {};
    % Purple continuation (UDUUDD shifted)
    \filldraw[cBlue] (0.5,0) -- (0.75,0.5) -- (1,0) -- (1.25,0.5) -- (1.5,1) -- (1.75,0.5) -- (2,0) -- cycle;
    \foreach \p in {(0.75,0.5),(1,0),(1.25,0.5),(1.5,1),(1.75,0.5),(2,0)} {
      \node[dot] at \p {};
    }
    \path (current bounding box.center) coordinate (-center);
  }
}

% 7. UUDDUUDD - starts with UU (intrinsic purple)
% Heights: 1,2,1,0,1,2,1,0
\tikzset{
  D4-UUDDUUDD/.pic={
    \DyckPic{cBlue}%
      {(0,0) -- (0.25,0.5) -- (0.5,1) -- (0.75,0.5) -- (1,0) -- (1.25,0.5) -- (1.5,1) -- (1.75,0.5) -- (2,0)}%
      {(0,0),(0.25,0.5),(0.5,1),(0.75,0.5),(1,0),(1.25,0.5),(1.5,1),(1.75,0.5),(2,0)}
  }
}

% 10. UUDUUDDD - starts with UU (intrinsic purple)
% Heights: 1,2,1,2,3,2,1,0
\tikzset{
  D4-UUDUUDDD/.pic={
    \DyckPic{cYellow}%
      {(0,0) -- (0.25,0.5) -- (0.5,1) -- (0.75,0.5) -- (1,1) -- (1.25,1.5) -- (1.5,1) -- (1.75,0.5) -- (2,0)}%
      {(0,0),(0.25,0.5),(0.5,1),(0.75,0.5),(1,1),(1.25,1.5),(1.5,1),(1.75,0.5),(2,0)}
  }
}

% --- BLUE FIBER (maps to UUDDUD) ---

% 4. UDUUDDUD - UD + UUDDUD
\tikzset{
  D4-UDUUDDUD/.pic={
    % Black triangle for UD prefix
    \filldraw[cBlack] (0,0) -- (0.25,0.5) -- (0.5,0) -- cycle;
    \node[dot] at (0,0) {};
    \node[dot] at (0.25,0.5) {};
    \node[dot] at (0.5,0) {};
    % Blue continuation (UUDDUD shifted)
    \filldraw[cPurple] (0.5,0) -- (0.75,0.5) -- (1,1) -- (1.25,0.5) -- (1.5,0) -- (1.75,0.5) -- (2,0) -- cycle;
    \foreach \p in {(0.75,0.5),(1,1),(1.25,0.5),(1.5,0),(1.75,0.5),(2,0)} {
      \node[dot] at \p {};
    }
    \path (current bounding box.center) coordinate (-center);
  }
}

% 11. UUUDDDUD - starts with UU (intrinsic blue)
% Heights: 1,2,3,2,1,0,1,0
\tikzset{
  D4-UUUDDDUD/.pic={
    \DyckPic{cPurple}%
      {(0,0) -- (0.25,0.5) -- (0.5,1) -- (0.75,1.5) -- (1,1) -- (1.25,0.5) -- (1.5,0) -- (1.75,0.5) -- (2,0)}%
      {(0,0),(0.25,0.5),(0.5,1),(0.75,1.5),(1,1),(1.25,0.5),(1.5,0),(1.75,0.5),(2,0)}
  }
}

% 12. UUUDDUDD - starts with UU (intrinsic blue)
% Heights: 1,2,3,2,1,2,1,0
\tikzset{
  D4-UUUDDUDD/.pic={
    \DyckPic{cBrown}%
      {(0,0) -- (0.25,0.5) -- (0.5,1) -- (0.75,1.5) -- (1,1) -- (1.25,0.5) -- (1.5,1) -- (1.75,0.5) -- (2,0)}%
      {(0,0),(0.25,0.5),(0.5,1),(0.75,1.5),(1,1),(1.25,0.5),(1.5,1),(1.75,0.5),(2,0)}
  }
}

% --- BROWN FIBER (maps to UUDUDD) ---

% 3. UDUUDUDD - UD + UUDUDD
\tikzset{
  D4-UDUUDUDD/.pic={
    % Black triangle for UD prefix
    \filldraw[cBlack] (0,0) -- (0.25,0.5) -- (0.5,0) -- cycle;
    \node[dot] at (0,0) {};
    \node[dot] at (0.25,0.5) {};
    \node[dot] at (0.5,0) {};
    % Brown continuation (UUDUDD shifted)
    \filldraw[cBrown] (0.5,0) -- (0.75,0.5) -- (1,1) -- (1.25,0.5) -- (1.5,1) -- (1.75,0.5) -- (2,0) -- cycle;
    \foreach \p in {(0.75,0.5),(1,1),(1.25,0.5),(1.5,1),(1.75,0.5),(2,0)} {
      \node[dot] at \p {};
    }
    \path (current bounding box.center) coordinate (-center);
  }
}

% 13. UUUDUDDD - starts with UU (intrinsic brown)
% Heights: 1,2,3,2,3,2,1,0
\tikzset{
  D4-UUUDUDDD/.pic={
    \DyckPic{cYellow}%
      {(0,0) -- (0.25,0.5) -- (0.5,1) -- (0.75,1.5) -- (1,1) -- (1.25,1.5) -- (1.5,1) -- (1.75,0.5) -- (2,0)}%
      {(0,0),(0.25,0.5),(0.5,1),(0.75,1.5),(1,1),(1.25,1.5),(1.5,1),(1.75,0.5),(2,0)}
  }
}

% --- YELLOW FIBER (maps to UUUDDD) ---

% 5. UDUUUDDD - UD + UUUDDD
\tikzset{
  D4-UDUUUDDD/.pic={
    % Black triangle for UD prefix
    \filldraw[cBlack] (0,0) -- (0.25,0.5) -- (0.5,0) -- cycle;
    \node[dot] at (0,0) {};
    \node[dot] at (0.25,0.5) {};
    \node[dot] at (0.5,0) {};
    % Yellow continuation (UUUDDD shifted)
    \filldraw[cYellow] (0.5,0) -- (0.75,0.5) -- (1,1) -- (1.25,1.5) -- (1.5,1) -- (1.75,0.5) -- (2,0) -- cycle;
    \foreach \p in {(0.75,0.5),(1,1),(1.25,1.5),(1.5,1),(1.75,0.5),(2,0)} {
      \node[dot] at \p {};
    }
    \path (current bounding box.center) coordinate (-center);
  }
}

% 14. UUUUDDDD - full pyramid (intrinsic yellow)
% Heights: 1,2,3,4,3,2,1,0
\tikzset{
  D4-UUUUDDDD/.pic={
    \DyckPic{cYellow}%
      {(0,0) -- (0.25,0.5) -- (0.5,1) -- (0.75,1.5) -- (1,2) -- (1.25,1.5) -- (1.5,1) -- (1.75,0.5) -- (2,0)}%
      {(0,0),(0.25,0.5),(0.5,1),(0.75,1.5),(1,2),(1.25,1.5),(1.5,1),(1.75,0.5),(2,0)}
  }
}
% ------------------------------------------------
% 2. Place the D3 lattice and arrows
% ------------------------------------------------

% Background orange region for D3, drawn first
\begin{scope}[on background layer]
  % \fill[cOrange] (-1.2,-0.3) -- (-1.7,1.5) -- (-1.7,3.4) --
  %                 (0,4.8) -- (2.0,2.2) -- (0.4,-0.3) -- cycle;
  \fill[cOrange] (0.45, -0.8) -- (-1.5, 0.8) -- (-1.5, 2.6) --
                  (0.5, 4.5) -- (2.3, 0.8) -- (0.9,-0.3) -- cycle;
\end{scope}

% Nodes (each "pic" lives around its own local origin)
\path (0,4.5) pic (D3_top)     {Dthree-UUUDDD}; %yellow
\path (-1.8,2.6) pic (D3_midL) {Dthree-UUDUDD}; %brown
\path ( 1.8,0.8) pic (D3_midR) {Dthree-UUDDUD}; %blue
\path (-1.8,0.8) pic (D3_lowL) {Dthree-UDUUDD}; %purple
\path (0,-0.8)  pic (D3_bot)   {Dthree-UDUDUD}; %green

% Hasse arrows (Tamari pentagon)
\draw[->] (0.119, -0.53)  -- (-1.5, 0.8);
\draw[->] (-1.5, 1.2) -- (-1.5, 2.6);
\draw[->] (-0.67, 3.385) -- (0.5, 4.5);
\draw[->] (0.8, -0.4)  -- (2.3, 0.8);
\draw[->] (2.0725, 1.257) -- (0.5, 4.5);

\node[below=3pt] at (D3_bot-center) {$D_3$};

\begin{scope}[on background layer]
  % \fill[cOrange] (-1.2,-0.3) -- (-1.7,1.5) -- (-1.7,3.4) --
  %                 (0,4.8) -- (2.0,2.2) -- (0.4,-0.3) -- cycle;
  \fill[cOrange] (-6.45, -0.8) -- (-7.5, 0.8) -- (-7.5, 2.6) --
                  (-6, 4.5) -- (-4.5, 0.8)  -- cycle;
\end{scope}

\path (-7,4.5) pic (D3_top)      {D4-UDUUUDDD}; %yellow
\path (-8.8,2.6) pic (D3_midL)   {D4-UDUUDUDD}; %brown
\path ( -5.8,0.8) pic (D3_midR)  {D4-UDUDUUDD}; %blue
\path (-8.8,0.8) pic (D3_lowL)   {D4-UDUUDDUD}; %purple
\path (-7,-0.8)  pic (D3_bot)    {D4-UDUDUDUD}; %green
\path (-12,0.8)  pic (D3_bot)    {D4-UUDDUDUD}; %green
\path (-13,2.4)  pic (D3_bot)    {D4-UUDUDDUD}; %green
\path (-13,4.0)  pic (D3_bot)    {D4-UUUDDDUD}; %green
\path (-13,6.0)  pic (D3_bot)    {D4-UUUDDUDD}; %green
\path (-12.75,8.0)  pic (D3_bot) {D4-UUUDUDDD}; %green
\path (-10,10.0)  pic (D3_bot)    {D4-UUUUDDDD}; %green
\path (-4.8,2.5)  pic (D3_bot)    {D4-UUDDUUDD}; %green
\path (-10,4.0)  pic (D3_bot)     {D4-UUDUDUDD}; %green
\path (-9,6.0)  pic (D3_bot)      {D4-UUDUUDDD}; %green

\draw[->] (-6.75, -0.3)  -- (-7.5, 0.8);
\draw[->] (-7.5, 1.2) -- (-7.5, 2.6);
\draw[->] (-7.1, 3.085) -- (-6, 4.5);
\draw[->] (-6.15, -0.55)  -- (-4.49, 0.8);
\draw[->] (-4.66, 1.223) -- (-6, 4.5);

\draw[->] (-6.75, -0.3) -- node[midway, above] {Loop} (-11,0.8);
\draw[->] (-11.75,1.25) -- (-11.75,2.4);
\draw[->] (-11.75,2.9) -- node[midway, right] {Loop} (-11.75,4);
\draw[->] (-11.75,4.5) -- (-11.75,6);
\draw[->] (-11.75,6.5) -- (-11.75,8);
\draw[->] (-11,8.5) -- node[midway, right] {Loop} (-9.25,10.0);
\draw[->] (-4.06, 1.223) -- node[midway, right] {Loop} (-4, 2.5);
% Change "--" to "to" and add the bend option
\draw[->] (-10.25,1.3) to[bend left=5] (-4, 2.5);
\draw[->] (-8.125, 1.2) to[bend right=10] node[midway, right] {Loop} (-11.75,4);
\draw[->] (-11.75,2.9) to[bend right=10] (-10,4.0);
\draw[->] (-4, 3) to[bend right=40] (-7.5,7.0);
\draw[->] (-10,4.0) to node[midway, right] {Loop} (-11.75,6);
\draw[->] (-7.5, 3.1) to[bend right=50] node[midway, right] {Loop} (-11.75,6);
\draw[->] (-9,5) -- (-9,6.0);
\draw[->] (-7.5,7.0) to[bend right=20] node[midway, right] {Loop} (-9.25,10.0);
\draw[->] (-6, 5.5) to[bend right=40] node[midway, right] {Loop} (-9.25,10.0);

% \node[below=4pt] at (-6.15, 2.4) {$D_4$};
\node[below=4pt, align=center] at (-2.6, 5) {$\implies$ \\ \scriptsize epimorphism of graphs};
% \path (0,4.5) pic (D3_top)     {Dthree-UUUDDD}; %yellow
% \path (-1.8,2.6) pic (D3_midL) {Dthree-UUDUDD}; %brown
% \path ( 1.8,0.8) pic (D3_midR) {Dthree-UUDDUD}; %blue
% \path (-1.8,0.8) pic (D3_lowL) {Dthree-UDUUDD}; %purple
% \path (0,-0.8)  pic (D3_bot)   {Dthree-UDUDUD}; %green
% \path (0,4.5) pic (D3_top)     {Dthree-UUUDDD}; %yellow
% \path (-1.8,2.6) pic (D3_midL) {Dthree-UUDUDD}; %brown
% \path ( 1.8,0.8) pic (D3_midR) {Dthree-UUDDUD}; %blue
% \path (-1.8,0.8) pic (D3_lowL) {Dthree-UDUUDD}; %purple

\end{tikzpicture}

\item We notice from \cite{MR1644479} the lattices $F_{n}^{\infty}$ (and $F_{n}^{p}$ (for any positive integer $p$)) form projective Fraisse family of graphs with monotone bonding epimorphism once we define the bonding epimorphism same as in the above example, making the corresponding inverse limits Kelley as above. The pic is given below.

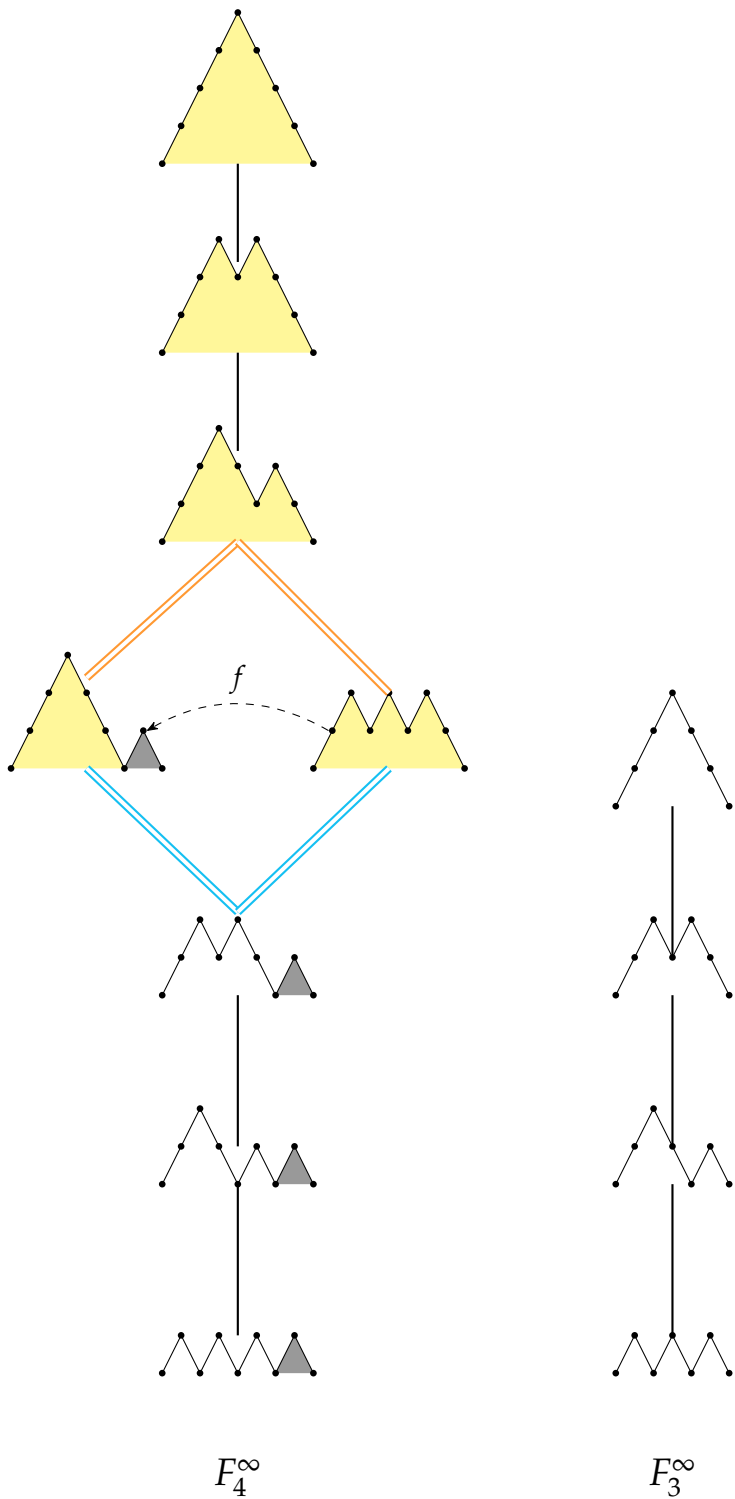
\tikzset{
  thick,
  >={Stealth[round]},
  dot/.style={circle,fill=black,inner sep=0pt,minimum size=2.5pt},
  % Colors
  cYellow/.style={fill=yellow!50},
  cGray/.style={fill=gray!80},
  % Edge styles
  edgeBlack/.style={draw=black, thick},
  edgeBlue/.style ={draw=cyan!70, double, double distance=1.5pt, thick},
  edgeOrange/.style={draw=orange!80, double, double distance=1.5pt, thick},
}

% =================================================================
% HELPER: Place Dots
% =================================================================
\newcommand{\PlaceDots}[1]{%
  \foreach \p in {#1} { \node[dot] at \p {}; }
}

% =================================================================
% PIC DEFINITIONS
% =================================================================

% --- RIGHT SIDE (F3) ---
\tikzset{
  F3-base/.pic={
    \draw (0,0) -- (0.25,0.5) -- (0.5,0) -- (0.75,0.5) -- (1,0) -- (1.25,0.5) -- (1.5,0);
    \PlaceDots{(0,0),(0.25,0.5),(0.5,0),(0.75,0.5),(1,0),(1.25,0.5),(1.5,0)}
    \coordinate (-top) at (0.75, 0.5); \coordinate (-bot) at (0.75, 0);
  }
}
\tikzset{
  F3-2/.pic={
    \draw (0,0) -- (0.25,0.5) -- (0.5,1) -- (0.75,0.5) -- (1,0) -- (1.25,0.5) -- (1.5,0);
    \PlaceDots{(0,0),(0.25,0.5),(0.5,1),(0.75,0.5),(1,0),(1.25,0.5),(1.5,0)}
    \coordinate (-top) at (0.75, 0.5); \coordinate (-bot) at (0.75, 0);
  }
}
\tikzset{
  F3-3/.pic={
    \draw (0,0) -- (0.25,0.5) -- (0.5,1) -- (0.75,0.5) -- (1,1) -- (1.25,0.5) -- (1.5,0);
    \PlaceDots{(0,0),(0.25,0.5),(0.5,1),(0.75,0.5),(1,1),(1.25,0.5),(1.5,0)}
    \coordinate (-top) at (0.75, 0.5); \coordinate (-bot) at (0.75, 0);
  }
}
\tikzset{
  F3-top/.pic={
    \draw (0,0) -- (0.25,0.5) -- (0.5,1) -- (0.75,1.5) -- (1,1) -- (1.25,0.5) -- (1.5,0);
    \PlaceDots{(0,0),(0.25,0.5),(0.5,1),(0.75,1.5),(1,1),(1.25,0.5),(1.5,0)}
    \coordinate (-top) at (0.75, 1.5); \coordinate (-bot) at (0.75, 0);
  }
}

% --- LEFT SIDE (F4) ---
\tikzset{
  F4-base/.pic={
    \fill[cGray] (1.5,0) -- (1.75,0.5) -- (2,0) -- cycle;
    \draw (0,0) -- (0.25,0.5) -- (0.5,0) -- (0.75,0.5) -- (1,0) -- (1.25,0.5) -- (1.5,0) -- (1.75,0.5) -- (2,0);
    \PlaceDots{(0,0),(0.25,0.5),(0.5,0),(0.75,0.5),(1,0),(1.25,0.5),(1.5,0),(1.75,0.5),(2,0)}
    \coordinate (-top) at (1, 0.5); \coordinate (-bot) at (1, 0);
  }
}
\tikzset{
  F4-2/.pic={
    \fill[cGray] (1.5,0) -- (1.75,0.5) -- (2,0) -- cycle;
    \draw (0,0) -- (0.25,0.5) -- (0.5,1) -- (0.75,0.5) -- (1,0) -- (1.25,0.5) -- (1.5,0) -- (1.75,0.5) -- (2,0);
    \PlaceDots{(0,0),(0.25,0.5),(0.5,1),(0.75,0.5),(1,0),(1.25,0.5),(1.5,0),(1.75,0.5),(2,0)}
    \coordinate (-top) at (1, 0.5); \coordinate (-bot) at (1, 0);
  }
}
\tikzset{
  F4-3/.pic={
    \fill[cGray] (1.5,0) -- (1.75,0.5) -- (2,0) -- cycle;
    \draw (0,0) -- (0.25,0.5) -- (0.5,1) -- (0.75,0.5) -- (1,1) -- (1.25,0.5) -- (1.5,0) -- (1.75,0.5) -- (2,0);
    \PlaceDots{(0,0),(0.25,0.5),(0.5,1),(0.75,0.5),(1,1),(1.25,0.5),(1.5,0),(1.75,0.5),(2,0)}
    \coordinate (-top) at (1, 1.1); 
    \coordinate (-bot) at (1, 0);
  }
}

% --- THE SPLIT ---
\tikzset{
  F4-SplitL/.pic={
    \fill[cYellow] (0,0) -- (0.25,0.5) -- (0.5,1) -- (0.75,1.5) -- (1,1) -- (1.25,0.5) -- (1.5,0) -- cycle;
    \fill[cGray] (1.5,0) -- (1.75,0.5) -- (2,0) -- cycle;
    \draw (0,0) -- (0.25,0.5) -- (0.5,1) -- (0.75,1.5) -- (1,1) -- (1.25,0.5) -- (1.5,0) -- (1.75,0.5) -- (2,0);
    \PlaceDots{(0,0),(0.25,0.5),(0.5,1),(0.75,1.5),(1,1),(1.25,0.5),(1.5,0),(1.75,0.5),(2,0)}
    \coordinate (-top) at (1, 1.2); \coordinate (-bot) at (1, 0);
    \coordinate (-linkR) at (1.8, 0.5);
  }
}
\tikzset{
  F4-SplitR/.pic={
    \fill[cYellow] (0,0) -- (0.25,0.5) -- (0.5,1) -- (0.75,0.5) -- (1,1) -- (1.25,0.5) -- (1.5,1) -- (1.75,0.5) -- (2,0) -- cycle;
    \draw (0,0) -- (0.25,0.5) -- (0.5,1) -- (0.75,0.5) -- (1,1) -- (1.25,0.5) -- (1.5,1) -- (1.75,0.5) -- (2,0);
    \PlaceDots{(0,0),(0.25,0.5),(0.5,1),(0.75,0.5),(1,1),(1.25,0.5),(1.5,1),(1.75,0.5),(2,0)}
    \coordinate (-top) at (1, 1); \coordinate (-bot) at (1, 0);
    \coordinate (-linkL) at (0.2, 0.5);
  }
}

% --- TOP STACK ---
\tikzset{
  F4-Top1/.pic={
    \fill[cYellow] (0,0) -- (0.25,0.5) -- (0.5,1) -- (0.75,1.5) -- (1,1) -- (1.25,0.5) -- (1.5,1) -- (1.75,0.5) -- (2,0) -- cycle;
    \draw (0,0) -- (0.25,0.5) -- (0.5,1) -- (0.75,1.5) -- (1,1) -- (1.25,0.5) -- (1.5,1) -- (1.75,0.5) -- (2,0);
    \PlaceDots{(0,0),(0.25,0.5),(0.5,1),(0.75,1.5),(1,1),(1.25,0.5),(1.5,1),(1.75,0.5),(2,0)}
    \coordinate (-top) at (1, 1.2); \coordinate (-bot) at (1, 0);
  }
}
\tikzset{
  F4-Top2/.pic={
    \fill[cYellow] (0,0) -- (0.25,0.5) -- (0.5,1) -- (0.75,1.5) -- (1,1) -- (1.25,1.5) -- (1.5,1) -- (1.75,0.5) -- (2,0) -- cycle;
    \draw (0,0) -- (0.25,0.5) -- (0.5,1) -- (0.75,1.5) -- (1,1) -- (1.25,1.5) -- (1.5,1) -- (1.75,0.5) -- (2,0);
    \PlaceDots{(0,0),(0.25,0.5),(0.5,1),(0.75,1.5),(1,1),(1.25,1.5),(1.5,1),(1.75,0.5),(2,0)}
    \coordinate (-top) at (1, 1.2); \coordinate (-bot) at (1, 0);
  }
}
\tikzset{
  F4-Top3/.pic={
    \fill[cYellow] (0,0) -- (1,2) -- (2,0) -- cycle;
    \draw (0,0) -- (0.25,0.5) -- (0.5,1) -- (0.75,1.5) -- (1,2) -- (1.25,1.5) -- (1.5,1) -- (1.75,0.5) -- (2,0);
    \PlaceDots{(0,0),(0.25,0.5),(0.5,1),(0.75,1.5),(1,2),(1.25,1.5),(1.5,1),(1.75,0.5),(2,0)}
    \coordinate (-top) at (1, 2); \coordinate (-bot) at (1, 0);
  }
}

% --- CORRECTED DUU SNIPPET (Down-Up-Up) ---
\tikzset{
  DUU-snippet/.pic={
    % D-U-U Path: Start high (0.5), go down to valley (0), go up (0.5), go up (1)
    \draw[thick, scale=0.7] (0,0.5) -- (0.25, 0) -- (0.5, 0.5) -- (0.75, 1);
    \foreach \p in {(0,0.5), (0.25,0), (0.5,0.5), (0.75,1)} 
      \node[dot, scale=0.7] at \p {};
  }
}

% =================================================================
% MAIN PICTURE
% =================================================================
\begin{tikzpicture}

  % --- F4 (LEFT COLUMN) ---
  \path (0,0)   pic (L1) {F4-base};
  \path (0,2.5) pic (L2) {F4-2};
  \path (0,5)   pic (L3) {F4-3};
  
  \path (-2, 8)  pic (LSplitL) {F4-SplitL};
  \path ( 2, 8)  pic (LSplitR) {F4-SplitR};
  
  \path (0, 11)   pic (L4) {F4-Top1};
  \path (0, 13.5) pic (L5) {F4-Top2};
  \path (0, 16)   pic (L6) {F4-Top3};

  % Edges for F4 Stack
  \draw[edgeBlack] (L1-top) -- (L2-bot);
  \draw[edgeBlack] (L2-top) -- (L3-bot);
  \draw[edgeBlack] (L4-top) -- (L5-bot);
  \draw[edgeBlack] (L5-top) -- (L6-bot);
  
  % Diamond Edges 
  \draw[edgeBlue]   (L3-top) -- (LSplitL-bot);
  \draw[edgeBlue]   (L3-top) -- (LSplitR-bot);
  \draw[edgeOrange] (LSplitL-top) -- (L4-bot);
  \draw[edgeOrange] (LSplitR-top) -- (L4-bot);
  
  % Arrow
  \draw[->, dashed, bend right=30] (LSplitR-linkL) to node[above] {$f$} (LSplitL-linkR);

  % --- F3 (RIGHT COLUMN) ---
  \path (6,0)   pic (R1) {F3-base};
  \path (6,2.5) pic (R2) {F3-2};
  \path (6,5)   pic (R3) {F3-3};
  \path (6,7.5) pic (R4) {F3-top};
  
  % Edges for F3
  \draw[edgeBlack] (R1-top) -- (R2-bot);
  \draw[edgeBlack] (R2-top) -- (R3-bot);
  \draw[edgeBlack] (R3-top) -- (R4-bot);

  % Labels
  \node[below=1cm] at (L1-bot) {\Large $F_4^\infty$};
  \node[below=1cm] at (R1-bot) {\Large $F_3^\infty$};

  % CORRECTED "avoid DUU" node with aligned baseline
  % \node[anchor=south east, align=right] at (7, 10) {
  %   avoid \tikz[baseline=0.3em]{\path pic {DUU-snippet};}
  % };

\end{tikzpicture}

\item Lattice of flats of Uniform matroids (Ranked, non distributive): Calling the lattice of flats of an uniform matroid as $U_{k,n}$ we notice that, there is an epimorphism of graphs from $U_{k,n}$ onto $U_{k,n-1}$ as follows: We notice that $U_{k,n}$ contains a copy of $U_{k,n-1}$. We first map the copy in the domain on to the range identically. For the rest in the second level from the top, we define \\
$\phi(a_{k}a_{k+1}\cdots a_{n-1}n)=a_{k+1}a_{k+2}\cdots a_{n-1}(n-1)$ if no $a_{i}=n-1$, $\phi(a_{k}a_{k+1}\cdots a_{n-1}n)=a_{k+1}a_{k+2}\cdots a_{n-2}(n-2)(n-1)$ if w.l.g. $a_{n-1}=n-1$\\
 $\phi(a_{k}a_{k+1}\cdots a_{n-1}n)=a_{k+1}a_{k+2}\cdots a_{n-3}(n-3)(n-2)(n-1)$ if w.l.g. $a_{n-2}=n-2, a_{n-1}=(n-1)$\\ in the third level from the top, we define $\phi(a_{k}a_{k+1}\cdots a_{n-2}n)=a_{k+1}a_{k+2}\cdots a_{n-2}(n-1)$ if no $a_{i}=n-1$= $\phi(a_{k}a_{k+1}\cdots a_{n-2}n)=a_{k+1}a_{k+2}\cdots a_{n-3}(n-2)(n-1)$ if w.l.g. $a_{n-2}=n-1$and it is equal to $a_{k+1}a_{k+2}\cdots a_{n-4}(n-3)(n-2)(n-1)$ if w.l.g. $a_{n-2}=n-2, a_{n-1}=(n-1)$, and so on. Proceeding in this way, we reach to the single digit $n$ and we define $\phi(n)=n-1$. We further notice that, each of these lattices of flats of $U_{k,n}$ are chain permutational posets. Again, each bonding map being confluent, the corresponding inverse limit is Kelley. 

\begin {figure}
\centering
\caption {The lattice of flats of uniform matroid}
\label {fig:poset}
\begin{tikzpicture}[
    scale=0.7,
    y=3cm, 
    x=0.6cm,
    every node/.style={font=\tiny\sffamily},
    dot/.style={circle, fill, inner sep=1.5pt},
    myred/.style={red!80!black},
    myblue/.style={blue!80!black},
    myblack/.style={black}
]

% =======================================================
% LEVEL 4: Full Set (1234567)
% =======================================================
% Positioned at x=18 (center of the lattice) and y=4.5
\node[dot, myred, label=above:{1234567}] (top) at (18, 4.5) {};

% =======================================================
% LEVEL 3: Triplets (ijk)
% =======================================================
\newcounter{countThree}

\foreach \i in {1,...,5} {
    \pgfmathtruncatemacro{\jstart}{\i+1}
    \foreach \j in {\jstart,...,6} {
        \pgfmathtruncatemacro{\kstart}{\j+1}
        \foreach \k in {\kstart,...,7} {
            \stepcounter{countThree}
            
            % Color Logic
            \ifnum\k=7 \def\col{myred} \else \ifnum\k=6 \def\col{myblue} \else \def\col{myblack} \fi \fi
            
            \node[dot, \col, label=above:{\i\j\k}] (n\i\j\k) at (\thecountThree, 3) {};
        }
    }
}

% =======================================================
% LEVEL 2: Pairs (xy)
% =======================================================
\newcounter{countTwo}

\foreach \x in {1,...,6} {
    \pgfmathtruncatemacro{\ystart}{\x+1}
    \foreach \y in {\ystart,...,7} {
        \stepcounter{countTwo}
        
        % Spacing adjustment
        \pgfmathsetmacro{\xpos}{\thecountTwo * 1.6 + 1} 
        
        % Color Logic
        \ifnum\y=7 \def\col{myred} \else \ifnum\y=6 \def\col{myblue} \else \def\col{myblack} \fi \fi
        
        \node[dot, \col, label=below:{\x\y}] (n\x\y) at (\xpos, 1.5) {};
    }
}

% =======================================================
% LEVEL 1: Singletons (z)
% =======================================================
\foreach \z in {1,...,7} {
    \pgfmathsetmacro{\xpos}{\z * 5 - 2}
    \ifnum\z=7 \def\col{myred} \else \ifnum\z=6 \def\col{myblue} \else \def\col{myblack} \fi \fi
    \node[dot, \col, label=below:{\z}] (n\z) at (\xpos, 0) {};
}

% =======================================================
% LEVEL 0: Empty Set
% =======================================================
% Positioned at x=18 (center) and y=-1.5
\node[dot, myblack, label=below:{$\emptyset$}] (bottom) at (18, -1.5) {};

% =======================================================
% EDGE GENERATION
% =======================================================
\begin{scope}[on background layer]
    
    % --- Edges: Level 4 (Top) to Level 3 ---
    % Since the top node {1..7} contains 7, all these edges are Red.
    \foreach \i in {1,...,5} {
        \pgfmathtruncatemacro{\jstart}{\i+1}
        \foreach \j in {\jstart,...,6} {
            \pgfmathtruncatemacro{\kstart}{\j+1}
            \foreach \k in {\kstart,...,7} {
                \draw[myred, opacity=0.4] (top.center) -- (n\i\j\k.center);
            }
        }
    }

    % --- Edges: Level 3 to Level 2 ---
    \foreach \i in {1,...,5} {
        \pgfmathtruncatemacro{\jstart}{\i+1}
        \foreach \j in {\jstart,...,6} {
            \pgfmathtruncatemacro{\kstart}{\j+1}
            \foreach \k in {\kstart,...,7} {
                
                \ifnum\k=7 \def\edgecol{myred} \else \ifnum\k=6 \def\edgecol{myblue} \else \def\edgecol{myblack} \fi \fi
                
                \draw[\edgecol, opacity=0.4] (n\i\j\k.center) -- (n\i\j.center);
                \draw[\edgecol, opacity=0.4] (n\i\j\k.center) -- (n\i\k.center);
                \draw[\edgecol, opacity=0.4] (n\i\j\k.center) -- (n\j\k.center);
            }
        }
    }

    % --- Edges: Level 2 to Level 1 ---
    \foreach \x in {1,...,6} {
        \pgfmathtruncatemacro{\ystart}{\x+1}
        \foreach \y in {\ystart,...,7} {
            
            \ifnum\y=7 \def\edgecol{myred} \else \ifnum\y=6 \def\edgecol{myblue} \else \def\edgecol{myblack} \fi \fi
            
            \draw[\edgecol, opacity=0.4] (n\x\y.center) -- (n\x.center);
            \draw[\edgecol, opacity=0.4] (n\x\y.center) -- (n\y.center);
        }
    }

    % --- Edges: Level 1 to Level 0 (Bottom) ---
    \foreach \z in {1,...,7} {
        % Color matches the singleton node above it
        \ifnum\z=7 \def\edgecol{myred} \else \ifnum\z=6 \def\edgecol{myblue} \else \def\edgecol{myblack} \fi \fi
        \draw[\edgecol, opacity=0.4] (n\z.center) -- (bottom.center);
    }

\end{scope}

\end{tikzpicture}
\end{figure}
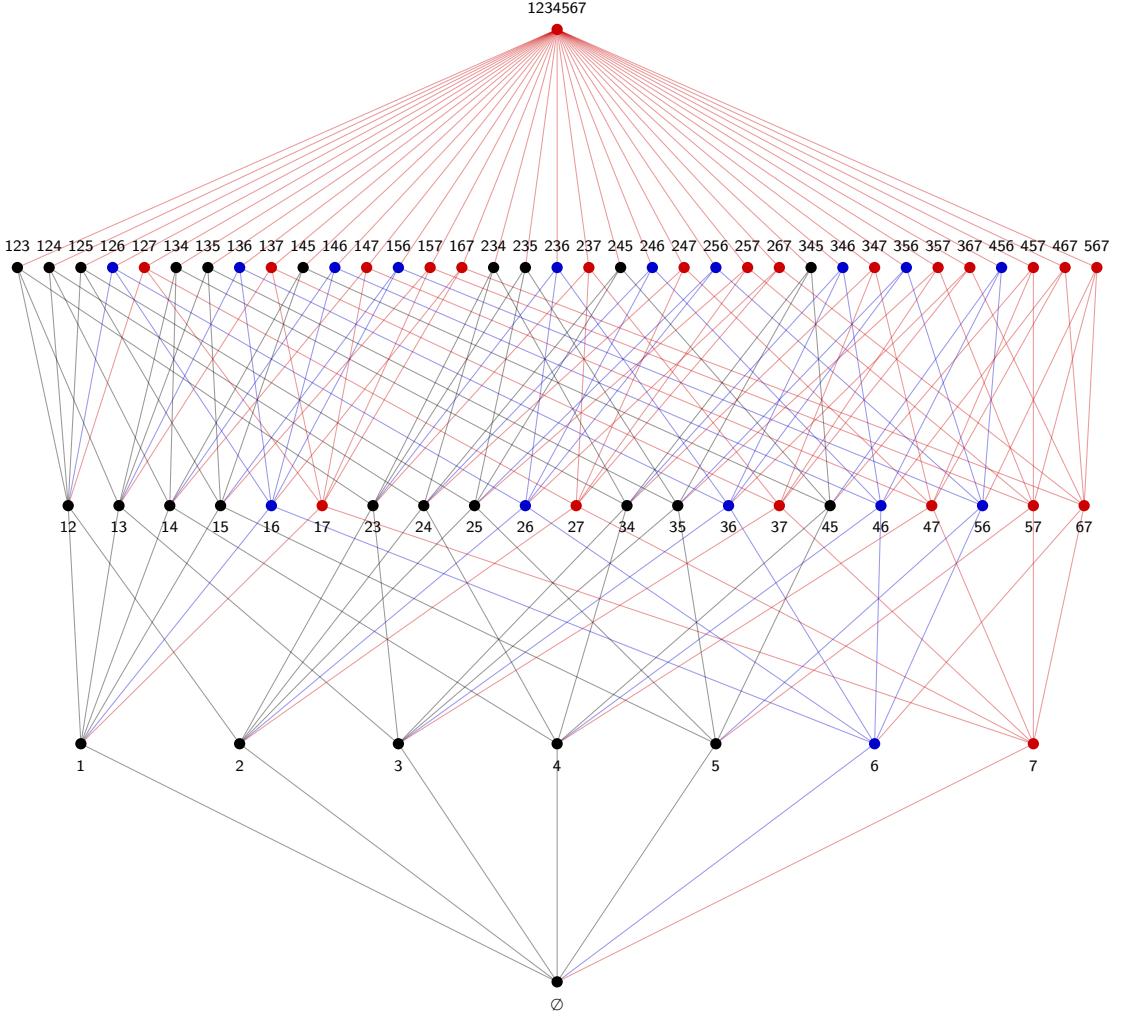

\item  Lattices of m dimer covers of snake graphs of continued fractions: Here we refer the snake graph of continuted fraction $G[a_{1}, a_{2}, \cdots a_{n}]$ as in as in $\cite{MR1644475}$ and $\cite{MR1644476}$ and we indicate the set of all $m$ dimer covers of them as in $\cite{MR1644475}$ and $\cite{MR1644476}$ as $G^{m}[a_{1}, a_{2}, \cdots a_{n}]$.  It is mentioned in section 2 of $\cite{MR1644476}$, that there is a bijection between the lattice of all $m$ dimer covers of the snake graph $G[a_{1}, a_{2}, \cdots a_{n}]$ and the $P-partition$ lattice of the shape  $G[a_{1}, a_{2}, \cdots a_{n}]$ with parts atmost $m$. Here we consider the scenario when $a_{1}=a_{2}=n$ creating a hook shape from the corresponding snake graph and we discuss in general the lattice $G^{m}[n,n]$. Note that, the Hasse diagram of $G^{m}[n,n]$ contains a copy of that of  $G^{m}[n-1,n-1]$ in it's top, where as it contains a copy of that of $G^{m-1}[n,n]$ in it's bottom. Let us call the copy of the smaller part as $Y$. So, we can define an epimorphism of graphs and posets from $G^{m}[n,n]$ on to $G^{m}[n-1,n-1]$ as follows: We map the copy of $G^{m}[n-1,n-1]$ in the top $Y$ identically. Then, we first notice that the Hasse diagram of $G^{m}[n,n]$ ($\forall m,n >1$) is symmetric w.r.t. the vertical line drawn through the top element of it. Suppose $x$ be a vertex in $G^{m}[n,n]$ out side $Y$. Then, if $x$ has same or smaller rank as the bottom most element of the copy of $Y$ , define $f(x)=y$, where $y$ is the bottom most element of $Y$. Now, if $x$ is outside $Y$, then whenever $x$ is to the right of the rightmost element in $Y$, say $x^{'}$ with the same rank, then define $f(x) =x^{'}$, if $x$ is to the left of the left most element in $Y$, say $x^{"}$ with the same rank, then define $f(x) =x^{"}$, if $x$ is somewhere in between the right most and leftmost elemnt in $Y$ with the same rank, then without loss of generality define $f(x) =x^{'}$. We can define an epimorphism of graphs and posets analogously from$G^{m}[n,n]$ on to $G^{m-1}[n,n]$. Thus, for any two pairs $(m,n) \geq (a,b)$, we can define the epimorphism of grapps and posets as  $\phi_{(a,b), (m,n)}$  from $G^{m}[n,n]$ on to $G^{a}[b,b]$ as the compositions of the intermidiate surjections making $(G^{m}[n,n], \phi_{(a,b),(m,n)})$ a projective Fraisse family of graphs preserving the partial order of the posets. Due to the vertical symmetry in the Hasse diagram, and the way the covering relations in the lattice are defined the corresponding bonding maps turns out to be confluent, making the graph of the projective limit Kelley as in $\cite{MR1644473}$. The necessary figure for this example and the next one are in the following page.
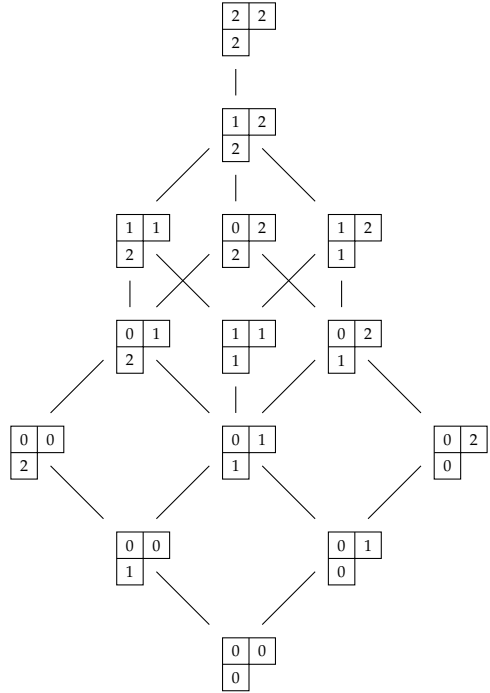
\begin {figure}[h]
\centering
\caption {The poset of $P$-partitions of $G[2,2]$ with parts at most $2$.}
\label {fig:poset}
\begin {tikzpicture}[scale=0.35]
    \newcommand{\drawsnake}{
        \draw (0,0) -- (0,1) -- (3,1) -- (3,0) -- cycle;
        \draw (1,0) -- (1,1);
        \draw (2,0) -- (2,1);
    }

    \newcommand{\drawdualsnake}{
        \draw (0,0) -- (0,2) -- (2,2) -- (2,1) -- (1,1) -- (1,0) -- cycle;
        \draw (0,1) -- (1,1) -- (1,2);
    }

    \newcommand{\drawposet}{
        \begin {scope}[shift={(-0.5,0)}]
        \draw[fill=black] (0,0) circle (0.05);
        \draw[fill=black] (1,1) circle (0.05);
        \draw[fill=black] (2,0) circle (0.05);
        \draw[fill=black] (0,1) circle (0.05);
        \draw[fill=black] (1,2) circle (0.05);
        \draw[fill=black] (2,1) circle (0.05);

        \draw (0,0) -- (1,1) -- (2,0);
        \draw (0,1) -- (1,2) -- (2,1);
        \draw (0,0) -- (0,1);
        \draw (1,1) -- (1,2);
        \draw (2,0) -- (2,1);
        \end {scope}
    }

    %%% Rank 0 %%%

    \drawdualsnake

    \draw (0.5,0.5) node {\tiny $0$};
    \draw (0.5,1.5) node {\tiny $0$};
    \draw (1.5,1.5) node {\tiny $0$};

    \draw (-0.5,2.5) -- (-2.5,4.5);
    \draw (1.5,2.5) -- (3.5,4.5);

    %%% Rank 1 %%%

    \begin {scope}[shift={(-4,4)}]
        \drawdualsnake

        \draw (0.5,0.5) node {\tiny $1$};
        \draw (0.5,1.5) node {\tiny $0$};
        \draw (1.5,1.5) node {\tiny $0$};

        \draw (-0.5,2.5) -- (-2.5,4.5);
        \draw (1.5,2.5) -- (3.5,4.5);
    \end {scope}

    \begin {scope}[shift={(4,4)}]
        \drawdualsnake

        \draw (0.5,0.5) node {\tiny $0$};
        \draw (0.5,1.5) node {\tiny $0$};
        \draw (1.5,1.5) node {\tiny $1$};

        \draw (-0.5,2.5) -- (-2.5,4.5);
        \draw (1.5,2.5) -- (3.5,4.5);
    \end {scope}

    %%% Rank 2 %%%

    \begin {scope}[shift={(-8,8)}]
        \drawdualsnake

        \draw (0.5,0.5) node {\tiny $2$};
        \draw (0.5,1.5) node {\tiny $0$};
        \draw (1.5,1.5) node {\tiny $0$};

        \draw (1.5,2.5) -- (3.5,4.5);
    \end {scope}

    \begin {scope}[shift={(0,8)}]
        \drawdualsnake

        \draw (0.5,0.5) node {\tiny $1$};
        \draw (0.5,1.5) node {\tiny $0$};
        \draw (1.5,1.5) node {\tiny $1$};

        \draw (0.5,2.5) -- (0.5,3.5);
        \draw (-0.5,2.5) -- (-2.5,4.5);
        \draw (1.5,2.5) -- (3.5,4.5);
    \end {scope}

    \begin {scope}[shift={(8,8)}]
        \drawdualsnake

        \draw (0.5,0.5) node {\tiny $0$};
        \draw (0.5,1.5) node {\tiny $0$};
        \draw (1.5,1.5) node {\tiny $2$};

        \draw (-0.5,2.5) -- (-2.5,4.5);
    \end {scope}

    %%% Rank 3 %%%

    \begin {scope}[shift={(-4,12)}]
        \drawdualsnake

        \draw (0.5,0.5) node {\tiny $2$};
        \draw (0.5,1.5) node {\tiny $0$};
        \draw (1.5,1.5) node {\tiny $1$};

        \draw (0.5,2.5) -- (0.5,3.5);
        \draw (1.5,2.5) -- (3.5,4.5);
    \end {scope}

    \begin {scope}[shift={(0,12)}]
        \drawdualsnake

        \draw (0.5,0.5) node {\tiny $1$};
        \draw (0.5,1.5) node {\tiny $1$};
        \draw (1.5,1.5) node {\tiny $1$};

        \draw (-0.5,2.5) -- (-2.5,4.5);
        \draw (1.5,2.5) -- (3.5,4.5);
    \end {scope}

    \begin {scope}[shift={(4,12)}]
        \drawdualsnake

        \draw (0.5,0.5) node {\tiny $1$};
        \draw (0.5,1.5) node {\tiny $0$};
        \draw (1.5,1.5) node {\tiny $2$};

        \draw (-0.5,2.5) -- (-2.5,4.5);
        \draw (0.5,2.5) -- (0.5,3.5);
    \end {scope}

    %%% Rank 4 %%%

    \begin {scope}[shift={(-4,16)}]
        \drawdualsnake

        \draw (0.5,0.5) node {\tiny $2$};
        \draw (0.5,1.5) node {\tiny $1$};
        \draw (1.5,1.5) node {\tiny $1$};

        \draw (1.5,2.5) -- (3.5,4.5);
    \end {scope}

    \begin {scope}[shift={(0,16)}]
        \drawdualsnake

        \draw (0.5,0.5) node {\tiny $2$};
        \draw (0.5,1.5) node {\tiny $0$};
        \draw (1.5,1.5) node {\tiny $2$};

        \draw (0.5,2.5) -- (0.5,3.5);
    \end {scope}

    \begin {scope}[shift={(4,16)}]
        \drawdualsnake

        \draw (0.5,0.5) node {\tiny $1$};
        \draw (0.5,1.5) node {\tiny $1$};
        \draw (1.5,1.5) node {\tiny $2$};

        \draw (-0.5,2.5) -- (-2.5,4.5);
    \end {scope}

    %%% Rank 5 %%%

    \begin {scope}[shift={(0,20)}]
        \drawdualsnake

        \draw (0.5,0.5) node {\tiny $2$};
        \draw (0.5,1.5) node {\tiny $1$};
        \draw (1.5,1.5) node {\tiny $2$};

        \draw (0.5,2.5) -- (0.5,3.5);
    \end {scope}

    %%% Rank 6 %%%

    \begin {scope}[shift={(0,24)}]
        \drawdualsnake

        \draw (0.5,0.5) node {\tiny $2$};
        \draw (0.5,1.5) node {\tiny $2$};
        \draw (1.5,1.5) node {\tiny $2$};
    \end {scope}
\end {tikzpicture}
\end {figure}

\item Lattices on Hessenberg functions: We notice that the poset of Hessenberg functions as in $\cite{MR1644490}$, $H_{n}$ contains a copy of that as $H_{n-1}$ to it's bottom right (there might be more copies). Hence, we can define an epimorphism from $H_{n}$ upon $H_{n-1}$ as follows: Lets map the copy of $H_{n-1}$ (call it as $H$) inside the domain on to the range identically. Then for any vertex $x$ outside $H$, if $x$ has a rank bigger or equal to that of $y$ (the top most element of $H$), then define $f(x)=y$, otherwise, if $x'$ is the leftmost element in $H$ with the same rank as that of $x$, define $f(x)=x^{'}$. The way the covering relations in the poset $H_{n}$ is defined, it makes the corresponding bonding maps confluent , giving the projective limit Kelley.

\begin{tikzpicture}[
    % Global Styles
    node distance=1.5cm,
    font=\small\ttfamily,
    >=Stealth,
    line width=0.6pt,
    xsep/.style={xshift=1.8cm},
    ysep/.style={yshift=1.2cm}
]

%% ---------------------------------------------------------
%% DIAGRAM (a): Hessenberg functions
%% ---------------------------------------------------------
\begin{scope}[local bounding box=diagramA]

    % --- Nodes (Positions remain the same) ---
    % Rank 0
    \node (1234) at (0,0) {1234};

    % Rank 1
    \node (1244) at (-1.8, 1.2) {1244};
    \node (1334) at (0, 1.2)    {1334};
    \node (2234) at (1.8, 1.2)  {2234};

    % Rank 2
    \node (1344) at (-1.8, 2.4) {1344};
    \node (2244) at (0, 2.4)    {2244};
    \node (2334) at (1.8, 2.4)  {2334};

    % Rank 3
    \node (1444) at (-1.8, 3.6) {1444};
    \node (2344) at (0, 3.6)    {2344};
    \node (3334) at (1.8, 3.6)  {3334};

    % Rank 4
    \node (2444) at (-0.9, 4.8) {2444};
    \node (3344) at (0.9, 4.8)  {3344};

    % Rank 5
    \node (3444) at (0, 6.0)    {3444};

    % Rank 6
    \node (4444) at (0, 7.2)    {4444};
    \node[right=0.5cm of 4444, font=\normalfont] {$h=4444$};

    % --- Edges (a) - REVERSED (High -> Low) ---

    % Rank 6 -> 5
    \draw[->] (4444) -- (3444);

    % Rank 5 -> 4
    \draw[->] (3444) -- (2444); % Reversed dotted
    \draw[->, dotted] (3444) -- (3344);

    % Rank 4 -> 3
    \draw[->] (2444) -- (1444);
    \draw[->, dotted] (2444) -- (2344);
    \draw[->, dotted] (3344) -- (2344);
    \draw[->, dotted] (3344) -- (3334);

    % Rank 3 -> 2
    \draw[->, dotted] (1444) -- (1344);
    \draw[->] (2344) -- (1344);
    \draw[->] (2344) -- (2244);
    \draw[->, dotted] (2344) -- (2334);
    \draw[->, dotted] (3334) -- (2334);

    % Rank 2 -> 1
    \draw[->] (1344) -- (1244); % Solid
    \draw[->, dotted] (1344) -- (1334);
    \draw[->, dotted] (2244) -- (1244);
    \draw[->, dotted] (2244) -- (2234);
    \draw[->] (2334) -- (1334);
    \draw[->, dotted] (2334) -- (2234);

    % Rank 1 -> 0
    \draw[->, dotted] (1244) -- (1234);
    \draw[->, dotted] (1334) -- (1234);
    \draw[->, dotted] (2234) -- (1234);

\end{scope}

\node[below=0.5cm of diagramA, font=\bfseries] {(a) Hessenberg functions.};

%% ---------------------------------------------------------
%% DIAGRAM (b): Degree tuples
%% ---------------------------------------------------------
\begin{scope}[xshift=6cm, local bounding box=diagramB]

    % --- Nodes (Positions remain the same) ---
    \node (1111) at (0,0) {1111};

    \node (2111) at (-1.8, 1.2) {2111};
    \node (1211) at (0, 1.2)    {1211};
    \node (1121) at (1.8, 1.2)  {1121};

    \node (2211) at (-1.8, 2.4) {2211};
    \node (2121) at (0, 2.4)    {2121};
    \node (1221) at (1.8, 2.4)  {1221};

    \node (3211) at (-1.8, 3.6) {3211};
    \node (2221) at (0, 3.6)    {2221};
    \node (1321) at (1.8, 3.6)  {1321};

    \node (3221) at (-0.9, 4.8) {3221};
    \node (2321) at (0.9, 4.8)  {2321};

    \node (3321) at (0, 6.0)    {3321};

    \node (4321) at (0, 7.2)    {4321};
    \node[right=0.5cm of 4321, font=\normalfont] {$\beta=4321$};

    % --- Edges (b) - REVERSED (High -> Low) ---

    % Explicit top edge to ensure visibility
    \draw[->] (4321) -- (3321);

    % Loop for the rest (reversed source/dest from previous)
    \foreach \high/\low in {
        3321/3221, 3321/2321,
        3221/3211, 3221/2221,
        2321/2221, 2321/1321,
        3211/2211,
        2221/2211, 2221/2121, 2221/1221,
        1321/2121, 1321/1221,
        2211/2111, 2211/1211,
        2121/2111, 2121/1121,
        1221/1211, 1221/1121,
        2111/1111, 1211/1111, 1121/1111} \draw[->] (\high) -- (\low);

\end{scope}

\node[below=0.5cm of diagramB, font=\bfseries] {(b) Degree tuples.};

\end{tikzpicture}

\item Lattices of non crossing partitions (no of blocks k=2): We notice that if $NC(n,k)$ is the set of all non crossing partitions of $[n]$ in to $k$ many blocks, then for $k =2$, always $NC(n,2)$ contains a copy of $NC(n-1,2)$ and we see a natural epimorphism from $NC(n,2)$ upon $NC(n-1,2)$, by utilising the technique of sequence of first occurrance $F(\rho)$ and it's complement $R(\rho)$ for any restricted grwoth functions $\rho$ as mentioned in section 7 of $\cite{MR1644478}$ to create the lattice structure on any $NC(n,k)$. By Lemma, 3.4 in \cite{MR1644473}we notice that here the bonding maps are monotone and by proposition 3.9 in \cite{MR1644473},if $N$ is the corresponding projective limit, then the corresponding projection of $N$ on to each finite factor $NC(n,2)$ is monotone. Each bonding map being monotone, by proposition 6.5 in \cite{MR1644473}, the corresponding projective limit $N$ is Kelley.\\

Note: As in case of the pattern avoidences, we can create graph sole-noid from those above examples also adjusting with specific subfamilies of the outer boundaries (and many alternatives).
 \begin{question}
\item 1. Are the projective Fraisse limit of the above lattices arcwise connected?
\item 2. Is this possible to create inverse system of monotone (or confluent) bonding maps over the family of Assecent lattice of $m$ dyck paths of length $n$ (and their mirrored versions) as in $\cite{MR1644477}$?
\item 3. Is this possible to create projective system of monotone (or confluent) bonding maps over the family of all $m$ dimer covers of all snake graphs of continued fraction $G[a_{1},a_{2},\cdots a_{n}]$ as in $\cite{MR1644475}, \cite{MR1644476}$?
\item 4. Is this possible to create projective system of monotone (or confluent) bonding maps over the family of all non crossing set partitions of $[n]$ with $k$ blocks varying both $n,k$ upto infinity (keeping $n>k$)?
\item 5.  Is this possible to create projective system of monotone (or confluent) bonding maps over the family of all lattice of flats of uniform matroids $U_{k,n}$varing  both $n,k$ upto infinity?

\end{question}
\end{enumerate}

Achknowledgement: Author thanks to Bruce Sagan for introducing the area of combinatorics. Thanks to Nicholas Ovenhouse for providing ideas and offering multiple helpful discussions, corrections in this topic. Thanks to Jean Luc Baril for many useful discussions, interest and several suggestions, corrections in the material. Thanks to Andrew Sack for idea and suggestions as well. Thanks to Mireille Bousquet-Melou, Richard Kenyon,  John Voight for suggestions and interest on this work. Thanks to Allen R. Williams for drawing all the pictures of the Hasse diagrams in Latex. Thanks to Robinson Czajkowski, David Gajewski for their many helful suggestions in latex.

Department of Mathematics and Statistics, University of Toledo, Ohio 43606, USA\\
*Correspondence:amrita.acharyya@utoledo.edu
\end{document}